\documentclass[12pt]{article}
\usepackage{full page, 
amssymb, amscd, graphicx, %showkeys
}
    \title{{\bf Logarithmic
tensor category theory, III: Intertwining maps and tensor
product bifunctors}}
    \author{Yi-Zhi Huang, James Lepowsky and Lin Zhang}
    \date{}

\begin{document}
    \bibliographystyle{alpha}
    \maketitle

    \newtheorem{rema}{Remark}[section]
    \newtheorem{propo}[rema]{Proposition}
    \newtheorem{theo}[rema]{Theorem}
   \newtheorem{defi}[rema]{Definition}
    \newtheorem{lemma}[rema]{Lemma}
    \newtheorem{corol}[rema]{Corollary}
     \newtheorem{exam}[rema]{Example}
\newtheorem{assum}[rema]{Assumption}
     \newtheorem{nota}[rema]{Notation}
        \newcommand{\ba}{\begin{array}}
        \newcommand{\ea}{\end{array}}
        \newcommand{\be}{\begin{equation}}
        \newcommand{\ee}{\end{equation}}
        \newcommand{\bea}{\begin{eqnarray}}
        \newcommand{\eea}{\end{eqnarray}}
        \newcommand{\nno}{\nonumber}
        \newcommand{\nn}{\nonumber\\}
        \newcommand{\lbar}{\bigg\vert}
        \newcommand{\p}{\partial}
        \newcommand{\dps}{\displaystyle}
        \newcommand{\bra}{\langle}
        \newcommand{\ket}{\rangle}
 \newcommand{\res}{\mbox{\rm Res}}
\newcommand{\wt}{\mbox{\rm wt}\;}
\newcommand{\swt}{\mbox{\scriptsize\rm wt}\;}
 \newcommand{\pf}{{\it Proof}\hspace{2ex}}
 \newcommand{\epf}{\hspace{2em}$\square$}
 \newcommand{\epfv}{\hspace{1em}$\square$\vspace{1em}}
        \newcommand{\ob}{{\rm ob}\,}
        \renewcommand{\hom}{{\rm Hom}}
\newcommand{\C}{\mathbb{C}}
\newcommand{\R}{\mathbb{R}}
\newcommand{\Z}{\mathbb{Z}}
\newcommand{\N}{\mathbb{N}}
\newcommand{\A}{\mathcal{A}}
\newcommand{\Y}{\mathcal{Y}}
\newcommand{\Arg}{\mbox{\rm Arg}\;}
\newcommand{\comp}{\mathrm{COMP}}
\newcommand{\lgr}{\mathrm{LGR}}

\newcommand{\dlt}[3]{#1 ^{-1}\delta \bigg( \frac{#2 #3 }{#1 }\bigg) }

\newcommand{\dlti}[3]{#1 \delta \bigg( \frac{#2 #3 }{#1 ^{-1}}\bigg) }

 \makeatletter
\newlength{\@pxlwd} \newlength{\@rulewd} \newlength{\@pxlht}
\catcode`.=\active \catcode`B=\active \catcode`:=\active \catcode`|=\active
\def\sprite#1(#2,#3)[#4,#5]{
   \edef\@sprbox{\expandafter\@cdr\string#1\@nil @box}
   \expandafter\newsavebox\csname\@sprbox\endcsname
   \edef#1{\expandafter\usebox\csname\@sprbox\endcsname}
   \expandafter\setbox\csname\@sprbox\endcsname =\hbox\bgroup
   \vbox\bgroup
  \catcode`.=\active\catcode`B=\active\catcode`:=\active\catcode`|=\active
      \@pxlwd=#4 \divide\@pxlwd by #3 \@rulewd=\@pxlwd
      \@pxlht=#5 \divide\@pxlht by #2
      \def .{\hskip \@pxlwd \ignorespaces}
      \def B{\@ifnextchar B{\advance\@rulewd by \@pxlwd}{\vrule
         height \@pxlht width \@rulewd depth 0 pt \@rulewd=\@pxlwd}}
      \def :{\hbox\bgroup\vrule height \@pxlht width 0pt depth
0pt\ignorespaces}
      \def |{\vrule height \@pxlht width 0pt depth 0pt\egroup
         \prevdepth= -1000 pt}
   }
\def\endsprite{\egroup\egroup}
\catcode`.=12 \catcode`B=11 \catcode`:=12 \catcode`|=12\relax
\makeatother

\def\hboxtr{\FormOfHboxtr} % Only necessary if
%\kern... is wanted
\sprite{\FormOfHboxtr}(25,25)[0.5 em, 1.2 ex] % Resolution ca. 200x340 dpi.

:BBBBBBBBBBBBBBBBBBBBBBBBB |
:BB......................B |
:B.B.....................B |
:B..B....................B |
:B...B...................B |
:B....B..................B |
:B.....B.................B |
:B......B................B |
:B.......B...............B |
:B........B..............B |
:B.........B.............B |
:B..........B............B |
:B...........B...........B |
:B............B..........B |
:B.............B.........B |
:B..............B........B |
:B...............B.......B |
:B................B......B |
:B.................B.....B |
:B..................B....B |
:B...................B...B |
:B....................B..B |
:B.....................B.B |
:B......................BB |
:BBBBBBBBBBBBBBBBBBBBBBBBB |

\endsprite

\def\shboxtr{\FormOfShboxtr} % Only necessary if
%\kern... is wanted
\sprite{\FormOfShboxtr}(25,25)[0.3 em, 0.72 ex] % Resolution ca. 200x340 dpi.

:BBBBBBBBBBBBBBBBBBBBBBBBB |
:BB......................B |
:B.B.....................B |
:B..B....................B |
:B...B...................B |
:B....B..................B |
:B.....B.................B |
:B......B................B |
:B.......B...............B |
:B........B..............B |
:B.........B.............B |
:B..........B............B |
:B...........B...........B |
:B............B..........B |
:B.............B.........B |
:B..............B........B |
:B...............B.......B |
:B................B......B |
:B.................B.....B |
:B..................B....B |
:B...................B...B |
:B....................B..B |
:B.....................B.B |
:B......................BB |
:BBBBBBBBBBBBBBBBBBBBBBBBB |

\endsprite

%\vspace{2em}

\begin{abstract}
This is the third part in a series of papers in which we introduce and
develop a natural, general tensor category theory for suitable module
categories for a vertex (operator) algebra.  In this paper (Part III),
we introduce and study intertwining maps and tensor product
bifunctors.
\end{abstract}

%\vspace{2em}

\tableofcontents
\vspace{2em}
%\noindent{\large \bf References}\hfill %290

In this paper, Part III of a series of eight papers on logarithmic
tensor category theory, we introduce and study intertwining maps and
tensor product bifunctors.  The sections, equations, theorems and so
on are numbered globally in the series of papers rather than within
each paper, so that for example equation (a.b) is the b-th labeled
equation in Section a, which is contained in the paper indicated as
follows: In Part I \cite{HLZ1}, which contains Sections 1 and 2, we
give a detailed overview of our theory, state our main results and
introduce the basic objects that we shall study in this work.  We
include a brief discussion of some of the recent applications of this
theory, and also a discussion of some recent literature.  In Part II
\cite{HLZ2}, which contains Section 3, we develop logarithmic formal
calculus and study logarithmic intertwining operators.  The present
paper, Part III, contains Section 4.  In Part IV \cite{HLZ4}, which
contains Sections 5 and 6, we give constructions of the $P(z)$- and
$Q(z)$-tensor product bifunctors using what we call ``compatibility
conditions'' and certain other conditions.  In Part V \cite{HLZ5},
which contains Sections 7 and 8, we study products and iterates of
intertwining maps and of logarithmic intertwining operators and we
begin the development of our analytic approach.  In Part VI
\cite{HLZ6}, which contains Sections 9 and 10, we construct the
appropriate natural associativity isomorphisms between triple tensor
product functors.  In Part VII \cite{HLZ7}, which contains Section 11,
we give sufficient conditions for the existence of the associativity
isomorphisms.  In Part VIII \cite{HLZ8}, which contains Section 12, we
construct braided tensor category structure.

\paragraph{Acknowledgments}
The authors gratefully
acknowledge partial support {}from NSF grants DMS-0070800 and
DMS-0401302.  Y.-Z.~H. is also grateful for partial support {}from NSF
grant PHY-0901237 and for the hospitality of Institut des Hautes 
\'{E}tudes Scientifiques in the fall of 2007.

\renewcommand{\theequation}{\thesection.\arabic{equation}}
\renewcommand{\therema}{\thesection.\arabic{rema}}
\setcounter{section}{3}
\setcounter{equation}{0}
\setcounter{rema}{0}

\section{$P(z)$- and $Q(z)$-intertwining maps and the
$P(z)$- and $Q(z)$-tensor product bifunctors}

We now generalize to the setting of the present work the notions of
$P(z)$- and $Q(z)$-tensor product of modules, for $z \in \C^{\times}$,
introduced in \cite{tensor1}, \cite{tensor2} and \cite{tensor3}. The symbols 
$P(z)$ and $Q(z)$ refer to moduli space elements described 
in Remarks \ref{P(z)geometry} and \ref{Q(z)geometry}, respectively. We
introduce the notions of $P(z)$- and $Q(z)$-intertwining map among
strongly $\tilde{A}$-graded generalized modules for a strongly
$A$-graded M\"{o}bius or conformal vertex algebra $V$ and establish
the relationship between such intertwining maps and grading-compatible
logarithmic intertwining operators.  We define the $P(z)$- and
$Q(z)$-tensor product bifunctors for pairs of strongly
$\tilde{A}$-graded generalized $V$-modules using these intertwining
maps and natural universal properties. As examples, for a strongly
$\tilde{A}$-graded generalized module $W$, we construct and describe
the $P(z)$-tensor products of $V$ and $W$ and also of $W$ and $V$; the
underlying strongly $\tilde{A}$-graded generalized modules of the
tensor product structures are $W$ itself, in both of these cases.  In
the case in which $V$ is a finitely reductive vertex operator algebra
(recall the Introduction), we construct and describe the $P(z)$- and
$Q(z)$-tensor products of arbitrary $V$-modules, and we use this
structure to motivate the construction of associativity isomorphisms
that we will carry out in later sections.  At the end of this section
we relate the $P(z)$- and $Q(z)$-tensor products.

We emphasize an important issue: Even though, as we have just
mentioned, we construct the $P(z)$- and $Q(z)$-tensor product
bifunctors in some cases, we do not give any {\it general}
construction of (models for) these bifunctors in this section.  But
for our deeper results, we will crucially need a suitable general
construction of these bifunctors, and indeed, for both $P(z)$ and
$Q(z)$, we will construct a useful, particular bifunctor (when it
exists) in Section 5.  We will use this construction in order to
construct the required natural associativity isomorphisms among triple
tensor products, leading to braided tensor category structure, under
suitable conditions.

In view of the results in Sections 2 and 3 involving contragredient
modules, it is natural for us to work in the strongly-graded setting
{}from now on:

\begin{assum}\label{assum}
Throughout this section and the remainder of this work, we shall
assume the following, unless other assumptions are explicitly made:
$A$ is an abelian group and $\tilde{A}$ is an abelian group containing
$A$ as a subgroup; $V$ is a strongly $A$-graded M\"{o}bius or
conformal vertex algebra; all $V$-modules and generalized $V$-modules
considered are strongly $\tilde{A}$-graded; and all intertwining
operators and logarithmic intertwining operators considered are
grading-compatible.  (Recall Definitions \ref{def:dgv}, \ref{def:dgw},
\ref{log:def} and \ref{gradingcompatintwop}.)
\end{assum}

We shall be working with full subcategories ${\cal C}$ of the category
${\cal M}_{sg}$ of strongly $\tilde{A}$-graded (ordinary) $V$-modules
or the category ${\cal GM}_{sg}$ of strongly $\tilde{A}$-graded
generalized $V$-modules (recall Notation \ref{MGM}).

In this section, $z$ will be a fixed nonzero complex number.

\subsection{$P(z)$-intertwining maps and the notion of $P(z)$-tensor product}

We first generalize the notion of $P(z)$-intertwining map given in
Section 4 of \cite{tensor1}; our $P(z)$-intertwining maps will
automatically be grading-compatible by definition.  We use the
notations given in Definition \ref{Wbardef}.  The main part of the
following definition, the Jacobi identity (\ref{im:def}), was
previewed in the Introduction (formula (\ref{im-jacobi})).  It should
be compared with the corresponding formula (\ref{intwmap}) in the Lie
algebra setting, and with the Jacobi identity (\ref{log:jacobi}) in
the definition of the notion of logarithmic intertwining operator;
note that the formal variable $x_2$ in that Jacobi identity is
specialized here to the nonzero complex number $z$.  Also, the
${\mathfrak s}{\mathfrak l}(2)$-bracket relations (\ref{im:Lj}) should
be compared with the corresponding relations (\ref{log:L(j)b}).  There
is no $L(-1)$-derivative formula for intertwining maps; as we shall
see, the $P(z)$-intertwining maps are obtained {}from logarithmic
intertwining operators by a process of specialization of the formal
variable to the complex variable $z$.

\begin{defi}\label{im:imdef}{\rm
Let $(W_1,Y_1)$, $(W_2,Y_2)$ and $(W_3,Y_3)$ be generalized
$V$-modules.  A {\it $P(z)$-intertwining map of type ${W_3\choose
W_1\,W_2}$} is a linear map
\begin{equation}\label{PzintwmapI}
I: W_1\otimes W_2 \to \overline{W}_3
\end{equation}
(recall {}from Definition \ref{Wbardef} that $\overline{W}_3$ is the
formal completion of $W_3$ with respect to the ${\mathbb C}$-grading)
such that the following conditions are satisfied: the {\it grading
compatibility condition}: for $\beta, \gamma\in \tilde{A}$ and
$w_{(1)}\in W_{1}^{(\beta)}$, $w_{(2)}\in W_{2}^{(\gamma)}$,
\begin{equation}\label{grad-comp}
I(w_{(1)}\otimes w_{(2)})\in \overline{W_{3}^{(\beta+\gamma)}};
\end{equation}
the
{\em lower truncation condition}: for any elements
$w_{(1)}\in W_1$, $w_{(2)}\in W_2$, and any $n\in {\mathbb C}$,
\begin{equation}\label{im:ltc}
\pi_{n-m}I(w_{(1)}\otimes w_{(2)})=0\;\;\mbox{ for }\;m\in {\mathbb N}
\;\mbox{ sufficiently large}
\end{equation}
(which follows {}from (\ref{grad-comp}), in view of the
grading restriction condition (\ref{set:dmltc}); recall the notation
$\pi_n$ {}from Definition \ref{Wbardef}); the {\em Jacobi identity}:
\begin{eqnarray}\label{im:def}
\lefteqn{x_0^{-1}\delta\bigg(\frac{ x_1-z}{x_0}\bigg)
Y_3(v, x_1)I(w_{(1)}\otimes w_{(2)})}\nno\\
&&=z^{-1}\delta\bigg(\frac{x_1-x_0}{z}\bigg)
I(Y_1(v, x_0)w_{(1)}\otimes w_{(2)})\nno\\
&&\hspace{2em}+x_0^{-1}\delta\bigg(\frac{z-x_1}{-x_0}\bigg)
I(w_{(1)}\otimes Y_2(v, x_1)w_{(2)})
\end{eqnarray}
for $v\in V$, $w_{(1)}\in W_1$ and $w_{(2)}\in W_2$ (note that all the
expressions in the right-hand side of (\ref{im:def}) are well defined,
and that the left-hand side of (\ref{im:def}) is meaningful because
any infinite linear combination of $v_n$ ($n\in{\mathbb Z}$) of the
form $\sum_{n<N}a_nv_n$ ($a_n\in {\mathbb C}$) acts in a well-defined
way on any $I(w_{(1)}\otimes w_{(2)})$, in view of (\ref{im:ltc}));
and the {\em ${\mathfrak s}{\mathfrak l}(2)$-bracket relations}: for any
$w_{(1)}\in W_1$ and $w_{(2)}\in W_2$,
\begin{equation}\label{im:Lj}
L(j)I(w_{(1)}\otimes w_{(2)})=I(w_{(1)}\otimes L(j)w_{(2)})+
\sum_{i=0}^{j+1}{j+1\choose i}z^iI((L(j-i)w_{(1)})\otimes w_{(2)})
\end{equation}
for $j=-1, 0$ and $1$ (note that if $V$ is in fact a conformal vertex
algebra, this follows automatically {}from (\ref{im:def}) by setting
$v=\omega$ and taking $\res_{x_0}\res_{x_1}x_1^{j+1}$). The vector
space of $P(z)$-intertwining maps of type ${W_3}\choose {W_1W_2}$ is
denoted by
\[
{\cal M}[P(z)]^{W_3}_{W_1W_2},
\]
or simply by
\[
{\cal M}^{W_3}_{W_1W_2}
\]
if there is no ambiguity. }
\end{defi}

\begin{rema}\label{P(z)geometry}
{\rm As we mentioned in the Introduction, $P(z)$ is the Riemann sphere
$\hat{\mathbb C}$ with one negatively oriented puncture at $\infty$
and two ordered positively oriented punctures at $z$ and $0$, with
local coordinates $1/w$, $w-z$ and $w$, respectively, vanishing at
these three punctures.  The geometry underlying the notion of
$P(z)$-intertwining map and the notions of $P(z)$-product and
$P(z)$-tensor product (see below) is determined by $P(z)$.}
\end{rema}

\begin{rema}{\rm
In the case of $\C$-graded ordinary modules for a vertex operator
algebra, where the grading restriction condition (\ref{Wn+k=0}) for a
module $W$ is replaced by the (more restrictive) condition
\begin{equation}
W_{(n)}=0 \;\; \mbox { for }\;n\in {\C}\;\mbox{ with sufficiently
negative real part}
\end{equation}
as in \cite{tensor1} (and where, in our context, the abelian groups
$A$ and $\tilde{A}$ are trivial), the notion of $P(z)$-intertwining
map above agrees with the earlier one introduced in \cite{tensor1}; in
this case, the conditions (\ref{grad-comp}) and (\ref{im:ltc}) are
automatic.}
\end{rema}

\begin{rema}\label{Pintwmaplowerbdd}
{\rm If $W_3$ in Definition \ref{im:imdef} is lower bounded, as in
Remark \ref{=0}, then (\ref{im:ltc}) can be strengthened to:
\begin{equation}\label{PpinI=0}
\pi_n I(w_{(1)}\otimes w_{(2)})=0\;\;\mbox{ for }
\;\Re{(n)}\;\mbox{ sufficiently negative}
\end{equation}
($n \in \C$).}
\end{rema}

\begin{rema}{\rm
As in Remark \ref{log:Lj2rema}, it is clear that the ${\mathfrak s}{\mathfrak
l}(2)$-bracket relations (\ref{im:Lj}) can equivalently be written as
\begin{eqnarray}\label{im:Lj2}
I(L(j)w_{(1)}\otimes w_{(2)})&=& \sum_{i=0}^{j+1}{j+1\choose
i}(-z)^i L(j-i)I(w_{(1)}\otimes w_{(2)})\nno\\
&&-\sum_{i=0}^{j+1}{j+1\choose i}(-z)^iI(w_{(1)}\otimes
L(j-i)w_{(2)})
\end{eqnarray}
for $w_{(1)}\in W_1$, $w_{(2)}\in W_2$ and $j=-1,0$ and $1$.
}
\end{rema}

Following \cite{tensor1} we will choose the branch of $\log z$ (and of
$\arg z$) such that
\begin{equation}\label{branch1}
0\leq \Im{(\log z)} = \arg z < 2 \pi
\end{equation}
(despite the fact that we happened to have used a different branch in
(\ref{log:br1}) in the proof of Theorem \ref{log:ids}), so that
\[
\log z = \log |z| + i \arg z. 
\]
We will also use the notation
\begin{equation}\label{branch2}
l_p(z)=\log z+2\pi ip, \ p\in {\mathbb Z},
\end{equation}
as in \cite{tensor1}, for arbitrary values of the $\log$ function. For
a formal expression $f(x)$ as in (\ref{log:f}), but involving only
nonnegative integral powers of $\log x$, and $\zeta\in {\mathbb C}$,
whenever
\begin{equation}\label{log:fsub}
f(x)\lbar_{x^n=e^{\zeta n},\;(\log x)^m=\zeta^m,\;n\in{\mathbb C},\;
m\in \N}
\end{equation}
exists algebraically, we will write (\ref{log:fsub}) simply as
$f(x)\lbar_{x=e^{\zeta}}$ or $f(e^\zeta)$, and we will call this
``substituting $e^\zeta$ for $x$ in $f(x)$,'' even though, in general,
it depends on $\zeta$, not just on $e^\zeta$.  (See also
(\ref{log:subs}).) In addition, for a fixed integer $p$, we will
sometimes write
\begin{equation}\label{im:f(z)}
f(x)\lbar_{x=z}\;\;\mbox{or}\;\;f(z)
\end{equation}
instead of $f(x)\lbar_{x=e^{l_p(z)}}$ or $f(e^{l_p(z)})$.  We will
sometimes say that ``$f(e^{\zeta})$ exists'' or that ``$f(z)$
exists.''

\begin{rema}{\rm
A very important example of an $f(z)$ existing in this sense occurs
when
\[
f(x)={\cal Y}(w_{(1)},x)w_{(2)} \;\; (\in W_3[\log x]\{x\})
\]
for $w_{(1)}\in W_1$, $w_{(2)}\in W_2$ and a logarithmic intertwining
operator ${\cal Y}$ of type ${W_3\choose W_1\,W_2}$, in the notation
of Definition \ref{log:def}; note that (\ref{log:fsub}) exists (as an
element of $\overline{W}_3$) in this case because of Proposition
\ref{log:logwt}(b).  Note also that in particular, ${\cal
Y}(w_{(1)},e^{\zeta})$ (or ${\cal Y}(w_{(1)},z)$) exists as a linear
map {}from $W_2$ to $\overline{W}_3$, and that ${\cal
Y}(\cdot,z)\cdot$ exists as a linear map
\begin{eqnarray}
W_1\otimes W_2 &\rightarrow & \overline{W}_3 \nonumber\\ 
w_{(1)}\otimes w_{(2)} &\mapsto & {\cal Y}(w_{(1)},z)w_{(2)}.
\end{eqnarray}
}
\end{rema}

Now we use these considerations to construct correspondences between
(grading-compatible) logarithmic intertwining operators and
$P(z)$-intertwining maps.  Fix an integer $p$. Let ${\cal Y}$ be a
logarithmic intertwining operator of type ${W_3\choose
W_1\,W_2}$.  Then we have a linear map
\begin{equation}
I_{{\cal Y},p}: W_1\otimes W_2\to \overline{W}_3
\end{equation}
defined by
\begin{equation}\label{log:IYp}
I_{{\cal Y},p}(w_{(1)}\otimes w_{(2)})={\cal
Y}(w_{(1)},e^{l_p(z)})w_{(2)}
\end{equation}
for all $w_{(1)}\in W_1$ and $w_{(2)}\in W_2$.  The
grading-compatibility condition (\ref{gradingcompatcondn}) yields the
grading-compatibility condition (\ref{grad-comp}) for $I_{{\cal
Y},p}$, and (\ref{im:ltc}) follows.  By substituting $e^{l_p(z)}$ for
$x_2$ in (\ref{log:jacobi}) and for $x$ in (\ref{log:L(j)b}), we see
that $I_{{\cal Y},p}$ satisfies the Jacobi identity (\ref{im:def}) and
the ${\mathfrak s}{\mathfrak l}(2)$-bracket relations (\ref{im:Lj}).
Hence $I_{{\cal Y},p}$ is a $P(z)$-intertwining map.  (Note that the
$L(-1)$-derivative property (\ref{log:L(-1)dev}) is not used here, so
that, for example, each ${\cal Y}^{(k)}$ in Remark \ref{Y(k)} produces
$P(z)$-intertwining maps in this way.  But the $L(-1)$-derivative
property is indeed needed for the recovery of ${\cal Y}$ from
$I_{{\cal Y},p}$, as we shall now see.)

On the other hand, we note that (\ref{log:p2}) (whose proof uses the
$L(-1)$-derivative property of ${\cal Y}$) is equivalent to
\begin{equation}\label{log:4.14}
\langle y^{L'(0)}w'_{(3)}, {\cal Y}(y^{-L(0)}w_{(1)},
x)y^{-L(0)}w_{(2)}\rangle_{W_3} =\langle w'_{(3)}, {\cal
Y}(w_{(1)}, xy)w_{(2)}\rangle_{W_3}
\end{equation}
for all $w_{(1)}\in W_1$, $w_{(2)}\in W_2$ and $w'_{(3)}\in W'_3$,
where we are using the pairing between the contragredient module
$W'_3$ and $W_3$ or $\overline{W}_3$ (recall Definition
\ref{defofWprime}, Theorem \ref{set:W'}, (\ref{L'(n)}),
(\ref{truncationforY'}) and (\ref{log:x^L(0)})).  Substituting
$e^{l_{p}(z)}$ for $x$ and then $e^{-l_{p}(z)}x$ for $y$, we obtain
\begin{eqnarray*}
&\langle y^{L'(0)}x^{L'(0)}w'_{(3)}, {\cal Y}
(y^{-L(0)}x^{-L(0)}w_{(1)}, e^{l_{p}(z)})
y^{-L(0)}x^{-L(0)}w_{(2)}\rangle_{W_3}\lbar_{y=e^{-l_{p}(z)}}&\\
&=\langle w'_{(3)}, {\cal Y}(w_{(1)}, x)w_{(2)}\rangle_{W_3},&
\end{eqnarray*}
or equivalently, using the notation (\ref{log:IYp}),
\begin{eqnarray*}
&\langle w'_{(3)}, y^{L(0)}x^{L(0)}I_{{\cal Y}, p}
(y^{-L(0)}x^{-L(0)}w_{(1)}\otimes
y^{-L(0)}x^{-L(0)}w_{(2)})\rangle_{W_3}\lbar_{y=e^{-l_{p}(z)}}&\\
&=\langle w'_{(3)}, {\cal Y}(w_{(1)}, x)w_{(2)}\rangle_{W_3}.&
\end{eqnarray*}
Thus we have recovered ${\cal Y}$ {}from $I_{{\cal Y},p}$ (with 
(\ref{log:L(-1)dev}) having been used in the proof).

This motivates the following definition: Given a $P(z)$-intertwining
map $I$ and an integer $p$, we define a linear map
\begin{equation}\label{YIp}
{\cal Y}_{I,p}:W_1\otimes W_2\to W_3[\log x]\{x\}
\end{equation}
by
\begin{eqnarray}\label{recover}
\lefteqn{{\cal Y}_{I,p}(w_{(1)}, x)w_{(2)}}\nno\\
&&=y^{L(0)}x^{L(0)}I(y^{-L(0)}x^{-L(0)}w_{(1)}\otimes
y^{-L(0)}x^{-L(0)}w_{(2)})\lbar_{y=e^{-l_{p}(z)}}
\end{eqnarray}
for any $w_{(1)}\in W_1$ and $w_{(2)}\in W_2$ (this is well defined
and indeed maps to $W_3[\log x]\{x\}$, in view of (\ref{log:x^L(0)})).
We will also use the notation ${w_{(1)}}_{n;k}^{I,p}w_{(2)}\in W_3$
defined by
\begin{equation}\label{wInkw}
{\cal Y}_{I,p}(w_{(1)}, x)w_{(2)}=\sum_{n\in{\mathbb
C}}\sum_{k\in {\mathbb N}} {w_{(1)}}_{n;k}^{I,p}w_{(2)}
x^{-n-1}(\log x)^k.
\end{equation}
Observe that since the operator $x^{\pm L(0)}$ always increases the
power of $x$ in an expression homogeneous of generalized weight $n$ by
$\pm n$, we see {}from (\ref{recover}) that
\begin{equation}\label{wt-cpnt-int-map}
{w_{(1)}}_{n;k}^{I,p}w_{(2)}\in(W_3)_{[n_1+n_2-n-1]}
\end{equation}
for $w_{(1)}\in (W_1)_{[n_1]}$ and $w_{(2)}\in (W_2)_{[n_2]}$.
Moreover, for $I=I_{{\cal Y},p}$, we have ${\cal Y}_{I,p}={\cal Y}$
(from the above), and for ${\cal Y}={\cal Y}_{I,p}$, we have $I_{{\cal
Y},p}=I$.

We can now prove the following proposition generalizing Proposition
12.2 in \cite{tensor3}.

\begin{propo}\label{im:correspond}
For $p\in {\mathbb Z}$, the correspondence
\[
{\cal Y}\mapsto I_{{\cal Y}, p}
\]
is a linear isomorphism {}from the space ${\cal
V}^{W_3}_{W_1W_2}$ of (grading-compatible) logarithmic intertwining
operators of type ${W_3\choose W_1\,W_2}$ to the space ${\cal
M}^{W_3}_{W_1W_2}$ of $P(z)$-intertwining maps of the same type. Its
inverse map is given by
\[
I\mapsto {\cal Y}_{I,p}.
\]
\end{propo}
\pf We need only show that for any $P(z)$-intertwining map $I$ of type
${W_3\choose W_1\,W_2}$, ${\cal Y}_{I,p}$ is a logarithmic
intertwining operator of the same type. The lower truncation condition 
(\ref{im:ltc}) implies that
the lower truncation condition (\ref{log:ltc}) for logarithmic
intertwining operator holds for ${\cal Y}_{I,p}$; for this,
(\ref{wt-cpnt-int-map}) can be used.  Let us now prove the
Jacobi identity for ${\cal Y}_{I,p}$.

Changing the formal variables $x_0$ and $x_1$ to
$x_0e^{l_p(z)}x_2^{-1}$ and $x_1e^{l_p(z)}x_2^{-1}$, respectively, in
the Jacobi identity (\ref{im:def}) for $I$, and then changing $v$ to
$y^{-L(0)}x_2^{-L(0)}v\lbar_{y=e^{-l_{p}(z)}}$ we obtain (noting that
at first, $e^{l_{p}(z)}$ could be written simply as $z$ because only
integral powers occur)
\begin{eqnarray*}
\lefteqn{x^{-1}_0\delta\left(\frac{x_1-x_2}{x_0}\right)
Y_{3}(y^{-L(0)}x_2^{-L(0)}v,x_1y^{-1}x_2^{-1}) I(w_{(1)}\otimes
w_{(2)})\lbar_{y=e^{-l_{p}(z)}}}\nno\\
&&=x_2^{-1}\delta\left(\frac{x_1-x_0}{x_2}\right)
I(Y_1(y^{-L(0)}x_2^{-L(0)}v, x_0y^{-1}x_2^{-1})w_{(1)}
\otimes w_{(2)})\lbar_{y=e^{-l_{p}(z)}}\nno\\
&&\quad +x_0^{-1}\delta\left(\frac{x_2-x_1}{-x_0}\right)
I(w_{(1)}\otimes Y_2(y^{-L(0)}x_2^{-L(0)}v, x_1
y^{-1}x_2^{-1})w_{(2)})\lbar_{y=e^{-l_{p}(z)}}.
\end{eqnarray*}
Using the formula
\[
Y_{3}(y^{-L(0)}x_2^{-L(0)}v,x_1y^{-1}x_2^{-1})=
y^{-L(0)}x_2^{-L(0)}Y_{3}(v,x_1)y^{L(0)}x_2^{L(0)},
\]
which holds on the generalized module $W_3$, by (\ref{log:p2}), and
the similar formulas for $Y_1$ and $Y_2$, we get
\begin{eqnarray*}
\lefteqn{x^{-1}_0\delta\left(\frac{x_1-x_2}{x_0}\right)
y^{-L(0)}x_2^{-L(0)}Y_{3}(v,x_1)y^{L(0)}x_2^{L(0)}I(w_{(1)}\otimes
w_{(2)})\lbar_{y=e^{-l_{p}(z)}}}\nno\\
&&=x_2^{-1}\delta\left(\frac{x_1-x_0}{x_2}\right)
I(y^{-L(0)}x_2^{-L(0)}Y_1(v, x_0)y^{L(0)}x_2^{L(0)}w_{(1)}
\otimes w_{(2)})\lbar_{y=e^{-l_{p}(z)}}\nno\\
&&\quad +x_0^{-1}\delta\left(\frac{x_2-x_1}{-x_0}\right)
I(w_{(1)}\otimes y^{-L(0)}x_2^{-L(0)}Y_2(v, x_1)
y^{L(0)}x_2^{L(0)}w_{(2)})\lbar_{y=e^{-l_{p}(z)}}.
\end{eqnarray*}
Replacing $w_{(1)}$ by
$y^{-L(0)}x_2^{-L(0)}w_{(1)}\lbar_{y=e^{-l_{p}(z)}}$ and $w_{(2)}$ by
$y^{-L(0)}x_2^{-L(0)}w_{(2)}\lbar_{y=e^{-l_{p}(z)}}$, and then
applying $y^{L(0)}x_2^{L(0)}\lbar_{y=e^{-l_{p}(z)}}$ to the whole
equation, we obtain
\begin{eqnarray*}
\lefteqn{x^{-1}_0\delta\left(\frac{x_1-x_2}{x_0}\right)
Y_{3}(v,x_1)y^{L(0)}x_2^{L(0)}\cdot}\nno\\
&&\hspace{2em}\cdot I(y^{-L(0)}x_2^{-L(0)}w_{(1)}\otimes
y^{-L(0)}x_2^{-L(0)}w_{(2)})\lbar_{y=e^{-l_{p}(z)}}\nno\\
&&=x_2^{-1}\delta\left(\frac{x_1-x_0}{x_2}\right)
y^{L(0)}x_2^{L(0)}\cdot\nno\\
&&\hspace{2em}\cdot I(y^{-L(0)}x_2^{-L(0)}Y_1(v, x_0)w_{(1)}\otimes
y^{-L(0)}x_2^{-L(0)}w_{(2)})\lbar_{y=e^{-l_{p}(z)}}\nno\\
&&\quad +x_0^{-1}\delta\left(\frac{x_2-x_1}{-x_0}\right)
y^{L(0)}x_2^{L(0)}\cdot\nno\\
&&\hspace{2em}\cdot I(y^{-L(0)}x_2^{-L(0)}w_{(1)}\otimes
y^{-L(0)}x_2^{-L(0)}Y_2(v,x_1) w_{(2)})\lbar_{y=e^{-l_{p}(z)}}.
\end{eqnarray*}
But using (\ref{recover}), we can write this as
\begin{eqnarray*}
\lefteqn{x^{-1}_0\delta\left(\frac{x_1-x_2}{x_0}\right)
Y_{3}(v, x_1){\cal Y}_{I, p}(w_{(1)}, x_2)w_{(2)}}\nno\\
&&=x_2^{-1}\delta\left(\frac{x_1-x_0}{x_2}\right) {\cal Y}_{I,
p}( Y_1(v, x_0) w_{(1)}, x_2)w_{(2)}\nno\\
&&\hspace{2em}+x_0^{-1}\delta\left(\frac{x_2-x_1}{-x_0}\right)
{\cal Y}_{I, p}(w_{(1)}, x_2) Y_2(v, x_1) w_{(2)}.
\end{eqnarray*}
That is, the Jacobi identity for ${\cal Y}_{I, p}$ holds.

Similar procedures show that the ${\mathfrak s}{\mathfrak l}(2)$-bracket
relations for $I$ imply the ${\mathfrak s}{\mathfrak l}(2)$-bracket
relations for ${\cal Y}_{I, p}$, as follows: Let $j$ be $-1$,
$0$ or $1$. By multiplying (\ref{im:Lj}) by $(yx)^j$ and using
(\ref{log:xLx^}) we obtain
\begin{eqnarray*}
\lefteqn{(yx)^{-L(0)}L(j)(yx)^{L(0)}I(w_{(1)}\otimes w_{(2)})}\nn
&&=I(w_{(1)}\otimes (yx)^{-L(0)}L(j)(yx)^{L(0)}w_{(2)})\nn
&&\quad +\sum_{i=0}^{j+1} {j+1\choose i}z^i(yx)^i
I(((yx)^{-L(0)}L(j-i)(yx)^{L(0)}w_{(1)})\otimes w_{(2)})
\end{eqnarray*}
Replacing $w_{(1)}$ by $(yx)^{-L(0)}w_{(1)}$ and $w_{(2)}$ by
$(yx)^{-L(0)}w_{(2)}$, and then applying $(yx)^{L(0)}$ to the whole
equation, we obtain
\begin{eqnarray*}
\lefteqn{L(j)(yx)^{L(0)}I((yx)^{-L(0)}w_{(1)}\otimes
(yx)^{-L(0)}w_{(2)})}\\
&&\hspace{3em}=(yx)^{L(0)}I((yx)^{-L(0)}w_{(1)}\otimes
(yx)^{-L(0)}L(j)w_{(2)})\\
&&\hspace{4em}+\sum_{i=0}^{j+1} {j+1\choose i}z^i(yx)^i
(yx)^{L(0)}I(((yx)^{-L(0)}L(j-i)w_{(1)})\otimes(yx)^{-L(0)}w_{(2)}).
\end{eqnarray*}
Evaluating at $y=e^{-l_p(z)}$ and using (\ref{recover}) we see that
this gives exactly the ${\mathfrak s}{\mathfrak l}(2)$-bracket relations
(\ref{log:L(j)b}) for ${\cal Y}_{I, p}$.

Finally, we prove the $L(-1)$-derivative property for ${\cal
Y}_{I,p}$. This follows {}from (\ref{recover}), (\ref{log:dx^}), and the
${\mathfrak s}{\mathfrak l}(2)$-bracket relation with $j=0$ for ${\cal
Y}_{I, p}$, namely,
\[
[L(0),{\cal Y}_{I, p}(w_{(1)},x)]={\cal Y}_{I, p}(L(0)w_{(1)},x)+
x{\cal Y}_{I, p}(L(-1)w_{(1)},x),
\]
as follows:
\begin{eqnarray*}
\lefteqn{\frac{d}{dx}{\cal Y}_{I, p}(w_{(1)},x)w_{(2)}}\\
&&=\frac{d}{dx}e^{-l_p(z)L(0)}x^{L(0)}I(e^{l_p(z)L(0)}x^{-L(0)}w_{(1)}
\otimes e^{l_p(z)L(0)}x^{-L(0)}w_{(2)})\\
&&=e^{-l_p(z)L(0)}x^{-1}x^{L(0)}L(0)I(e^{l_p(z)L(0)}x^{-L(0)}w_{(1)}
\otimes e^{l_p(z)L(0)}x^{-L(0)}w_{(2)})\\
&&\hspace{1em}-e^{-l_p(z)L(0)}x^{L(0)}I(e^{l_p(z)L(0)}x^{-1}x^{-L(0)}
L(0)w_{(1)}\otimes e^{l_p(z)L(0)}x^{-L(0)}w_{(2)})\\
&&\hspace{1em}-e^{-l_p(z)L(0)}x^{L(0)}I(e^{l_p(z)L(0)}x^{-L(0)}w_{(1)}
\otimes e^{l_p(z)L(0)}x^{-1}x^{-L(0)}L(0)w_{(2)})\\
&&=x^{-1}L(0){\cal Y}_{I, p}(w_{(1)},x)w_{(2)}-x^{-1}{\cal Y}_{I, p}
(w_{(1)},x)L(0)w_{(2)}\\
&&\hspace{1em}-x^{-1}{\cal Y}_{I, p}(L(0)w_{(1)},x)w_{(2)}\\
&&={\cal Y}_{I, p}(L(-1)w_{(1)},x)w_{(2)}. \hspace{16em}\square
\end{eqnarray*}

\vspace{.1em}

\begin{rema}\label{lowerbddcorrespondence}{\rm
{}From Remarks \ref{=0} and \ref{Pintwmaplowerbdd}, we note that if
$W_3$ is lower bounded, then the spaces of logarithmic intertwining
operators and of $P(z)$-intertwining maps in Proposition
\ref{im:correspond} satisfy the stronger conditions
(\ref{repartbounded}) and (\ref{PpinI=0}), respectively.}

\end{rema}

\begin{rema}\label{mod-sub}{\rm
Given a generalized $V$-module $(W, Y_{W})$, recall {}from Remark
\ref{str-graded-g-mod-as-l-int} that $Y_{W}$ is a logarithmic
intertwining operator of type ${W\choose VW}$ not involving $\log x$
and having only integral powers of $x$. Then the substitution
$x\mapsto z$ in (\ref{log:IYp}) is very simple; it is independent of
$p$ and $Y_{W}(\cdot, z)\cdot$ entails only the substitutions
$x^{n}\mapsto z^{n}$ for $n\in \Z$. As a special case, we can take
$(W, Y_{W})$ to be $(V, Y)$ itself.}
\end{rema}

\begin{rema}\label{Ypp'}{\rm
Let $I$ be a $P(z)$-intertwining map of type ${W_3\choose W_1\,W_2}$
and let $p,p' \in {\mathbb Z}$.  {}From (\ref{recover}), we see that the
logarithmic intertwining operators ${\cal Y}_{I, p}$ and ${\cal Y}_{I,
p'}$ of this same type differ as follows:
\begin{eqnarray}\label{YIp'YIp}
\lefteqn{{\cal Y}_{I,p'}(w_{(1)}, x)w_{(2)}}\nno\\
&&=e^{2\pi i(p-p')L(0)}{\cal Y}_{I,p}(e^{2\pi i(p'-p)L(0)}w_{(1)},
x)e^{2\pi i(p'-p)L(0)}w_{(2)}
\end{eqnarray}
for $w_{(1)}\in W_1$ and $w_{(2)}\in W_2$.  Using the notation in
Remark \ref{Ys1s2s3}, we thus have
\begin{eqnarray}
{\cal Y}_{I,p'}&=&({\cal Y}_{I,p})_{[p-p',p'-p,p'-p]}\nno\\
&=&{\cal Y}_{I,p}(\cdot,e^{2\pi i(p-p')} \cdot)\cdot .
\end{eqnarray}
}
\end{rema}

\begin{rema}\label{II1}{\rm
Let $I$ be a $P(z)$-intertwining map of type ${W_3}\choose
{W_1W_2}$. Then {}from the correspondence between $P(z)$-intertwining
maps and logarithmic intertwining operators in Proposition
\ref{im:correspond}, we see that for any nonzero complex number $z_1$,
the linear map $I_1$ defined by
\begin{equation}\label{log:zz_1}
I_1(w_{(1)}\otimes w_{(2)})=\sum_{n\in{\mathbb C}}\sum_{k\in {\mathbb N}}
{w_{(1)}}_{n;k}^{I,p}w_{(2)} e^{l_p(z_1)(-n-1)}(l_p(z_1))^k
\end{equation}
for $w_{(1)}\in W_1$ and $w_{(2)}\in W_2$ (recall (\ref{wInkw})) is a
$P(z_1)$-intertwining map of the same type.  In this sense,
${w_{(1)}}_{n;k}^{I,p}w_{(2)}$ is independent of $z$.  This justifies
writing $I(w_{(1)}\otimes w_{(2)})$ alternatively as
\begin{eqnarray}\label{imz}
I(w_{(1)},z)w_{(2)},
\end{eqnarray}
indicating that $z$ can be replaced by any nonzero complex number;
this notation was sometimes used in \cite{tensor4}, although we shall
generally not be using it in the present work.  However, for a general
intertwining map associated to a sphere with punctures not necessarily
of type $P(z)$, the corresponding element
${w_{(1)}}_{n;k}^{I,p}w_{(2)}$ will in general be different.  }
\end{rema}

We now proceed to the definition of the $P(z)$-tensor product.  As in
\cite{tensor1}, this will be a suitably universal ``$P(z)$-product.''
We generalize these notions {}from \cite{tensor1} using the notations
${\cal M}_{sg}$ and ${\cal GM}_{sg}$ (the categories of strongly
graded $V$-modules and generalized $V$-modules, respectively; recall
Notation \ref{MGM}) as follows:

\begin{defi}\label{pz-product}{\rm
Let ${\cal C}_1$ be either of the categories ${\cal M}_{sg}$ or ${\cal
GM}_{sg}$ (recall Notation \ref{MGM}).  For $W_1, W_2\in \ob{\cal
C}_1$, a {\em $P(z)$-product of $W_1$ and $W_2$} is an object
$(W_3,Y_3)$ of ${\cal C}_1$ equipped with a $P(z)$-intertwining map
$I_3$ of type ${W_3\choose W_1\,W_2}$. We denote it by $(W_3,Y_3;I_3)$
or simply by $(W_3;I_3)$. Let $(W_4,Y_4;I_4)$ be another
$P(z)$-product of $W_1$ and $W_2$. A {\em morphism} {}from
$(W_3,Y_3;I_3)$ to $(W_4,Y_4;I_4)$ is a module map $\eta$ {}from $W_3$
to $W_4$ such that the diagram
\begin{center}
\begin{picture}(100,60)
\put(-5,0){$\overline W_3$}
\put(13,4){\vector(1,0){104}}
\put(119,0){$\overline W_4$}
\put(41,50){$W_1\otimes W_2$}
\put(61,45){\vector(-3,-2){50}}
\put(68,45){\vector(3,-2){50}}
\put(65,8){$\bar\eta$}
\put(20,27){$I_3$}
\put(98,27){$I_4$}
\end{picture}
\end{center}
commutes, that is,
\begin{equation}
I_4=\bar\eta\circ I_3,
\end{equation}
where
\begin{equation}
\bar\eta : \overline{W}_3 \to \overline{W}_4
\end{equation}
is the natural extension of $\eta$.  (Note that $\bar\eta$
exists because $\eta$ preserves ${\mathbb C}$-gradings; we shall use
the notation $\bar\eta$ for any such map $\eta$.)
}
\end{defi}

\begin{rema}{\rm
In this setting, let $\eta$ be a morphism {}from $(W_3,Y_3;I_3)$ to
$(W_4,Y_4;I_4)$.  We know {}from (\ref{YIp})--(\ref{wInkw}) that for $p
\in {\mathbb Z}$, the coefficients ${w_{(1)}}_{n;k}^{I_{3},p}w_{(2)}$
and ${w_{(1)}}_{n;k}^{I_{4},p}w_{(2)}$ in the formal expansion
(\ref{wInkw}) of ${\cal Y}_{I_{3},p}(w_{(1)}, x)w_{(2)}$ and ${\cal
Y}_{I_{4},p}(w_{(1)}, x)w_{(2)}$, respectively, are determined by
$I_3$ and $I_4$, and that
\begin{equation}\label{etaw1w2}
\eta ({w_{(1)}}_{n;k}^{I_{3},p}w_{(2)}) =
{w_{(1)}}_{n;k}^{I_{4},p}w_{(2)},
\end{equation}
as we see by applying $\bar\eta$ to (\ref{recover}).
}
\end{rema}

The notion of $P(z)$-tensor product is now defined by means of a
universal property as follows:

\begin{defi}\label{pz-tp}{\rm
Let ${\cal C}$ be a full subcategory of either ${\cal M}_{sg}$ or
${\cal GM}_{sg}$.  For $W_1, W_2\in \ob{\cal C}$, a {\em $P(z)$-tensor
product of $W_1$ and $W_2$ in ${\cal C}$} is a $P(z)$-product $(W_0,
Y_0; I_0)$ with $W_0\in{\rm ob\,}{\cal C}$ such that for any
$P(z)$-product $(W,Y;I)$ with $W\in{\rm ob\,}{\cal C}$, there is a
unique morphism {}from $(W_0, Y_0; I_0)$ to $(W,Y;I)$.  Clearly, a
$P(z)$-tensor product of $W_1$ and $W_2$ in ${\cal C}$, if it exists,
is unique up to unique isomorphism. In this case we will denote it by
\[
(W_1\boxtimes_{P(z)} W_2, Y_{P(z)}; \boxtimes_{P(z)})
\]
and call the object
\[
(W_1\boxtimes_{P(z)} W_2, Y_{P(z)})
\]
the {\em $P(z)$-tensor product (generalized) module of $W_1$ and $W_2$ in
${\cal C}$}.  We will skip the phrase ``in ${\cal C}$'' if the
category ${\cal C}$ under consideration is clear in context.  }
\end{defi}

\begin{rema}{\rm
Consider the functor {}from ${\cal C}$ to the category ${\bf Set}$
defined by assigning to $W\in {\rm ob\,}{\cal C}$ the set ${\cal
M}^W_{W_1\,W_2}$ of all $P(z)$-intertwining maps of type ${W\choose
W_1\,W_2}$.  Then if the $P(z)$-tensor product of $W_1$ and $W_2$
exists, it is just the universal element for this functor, and this
functor is representable, represented by the $P(z)$-tensor product.
(Recall that given a functor $f$ {}from a category ${\cal K}$ to ${\bf
Set}$, a universal element for $f$, if it exists, is a pair $(X,x)$
where $X\in{\rm ob\,}{\cal K}$ and $x\in f(X)$ such that for any pair
$(Y,y)$ with $Y\in {\rm ob\,}{\cal K}$ and $y\in f(Y)$, there is a
unique morphism $\sigma: X\to Y$ such that $f(\sigma)(x)=y$; in this
case, $f$ is represented by $X$.)  }
\end{rema}

Definition \ref{pz-tp} and Proposition \ref{im:correspond} immediately
give the following result relating the module maps {}from a
$P(z)$-tensor product (generalized) module with the
$P(z)$-intertwining maps and the logarithmic intertwining operators:

\begin{propo}\label{pz-iso}
Suppose that $W_1\boxtimes_{P(z)}W_2$ exists. We have a natural
isomorphism
\begin{eqnarray}\label{isofromhomstointwmaps}
\hom_{V}(W_1\boxtimes_{P(z)}W_2, W_3)&\stackrel{\sim}{\to}&
{\cal M}^{W_3}_{W_1W_2}\nno\\
\eta&\mapsto& \overline{\eta}\circ
\boxtimes_{P(z)}
\end{eqnarray}
and for $p\in {\mathbb Z}$, a natural isomorphism
\begin{eqnarray}
\hom_{V}(W_1\boxtimes_{P(z)} W_2, W_3)&
\stackrel{\sim}{\rightarrow}& {\cal V}^{W_3}_{W_1W_2}\nno\\
\eta&\mapsto & {\cal Y}_{\eta, p}
\end{eqnarray}
where ${\cal Y}_{\eta, p}={\cal Y}_{I, p}$ with
$I=\overline{\eta}\circ \boxtimes_{P(z)}$.\epf
\end{propo}

Suppose that the $P(z)$-tensor product $(W_1\boxtimes_{P(z)} W_2,
Y_{P(z)}; \boxtimes_{P(z)})$ of $W_1$ and $W_2$ exists.  We will
sometimes denote the action of the canonical $P(z)$-intertwining map
\begin{equation}\label{actionofboxtensormap}
w_{(1)} \otimes w_{(2)}\mapsto\boxtimes_{P(z)}(w_{(1)} \otimes
w_{(2)})=\boxtimes_{P(z)}(w_{(1)}, z)w_{(2)} \in 
\overline{W_1\boxtimes_{P(z)}W_2}
\end{equation}
(recall (\ref{imz})) on elements simply by $w_{(1)}\boxtimes_{P(z)}
w_{(2)}$:
\begin{equation}\label{boxtensorofelements}
w_{(1)}\boxtimes_{P(z)} w_{(2)}=\boxtimes_{P(z)}(w_{(1)} \otimes
w_{(2)})=\boxtimes_{P(z)}(w_{(1)}, z)w_{(2)}.
\end{equation}

\begin{rema} {\rm
We emphasize that the element $w_{(1)}\boxtimes_{P(z)} w_{(2)}$
defined here is an element of the formal completion
$\overline{W_1\boxtimes_{P(z)}W_2}$, and {\it not} (in general) of the
module $W_1\boxtimes_{P(z)}W_2$ itself.  This is different {}from the
classical case for modules for a Lie algebra (recall Section
\ref{LA}), where the tensor product of elements of two modules is an
element of the tensor product module.}
\end{rema}

\begin{rema} {\rm
Note that under the natural isomorphism (\ref{isofromhomstointwmaps})
for the case $W_{3} = W_1\boxtimes_{P(z)}W_2$, the identity map {}from
$W_1\boxtimes_{P(z)}W_2$ to itself corresponds to the canonical
intertwining map $\boxtimes_{P(z)}$.  Furthermore, for $p \in {\mathbb
Z}$, the $P(z)$-tensor product of $W_1$ and $W_2$ gives rise to a
logarithmic intertwining operator ${\cal Y}_{\boxtimes_{P(z)}, p}$ of
type ${W_1\boxtimes_{P(z)}W_2\choose W_1\,W_2}$, according to formula
(\ref{recover}).  If $p$ is changed to $p' \in {\mathbb Z}$, this
logarithmic intertwining operator changes according to
(\ref{YIp'YIp}).  Note that the $P(z)$-intertwining map
$\boxtimes_{P(z)}$ is canonical and depends only on $z$, while a
corresponding logarithmic intertwining operator is not; it depends on
$p \in {\mathbb Z}$.
}
\end{rema}

\begin{rema} {\rm
Sometimes it will be convenient, as in the next proposition, to use
the particular isomorphism associated with $p=0$ (in Proposition
\ref{im:correspond}) between the spaces of $P(z)$-intertwining maps
and of logarithmic intertwining operators of the same type.  In this
case, we shall sometimes simplify the notation by dropping the $p$
($=0$) in the notation ${w_{(1)}}_{n;k}^{I,0}w_{(2)}$ (recall
(\ref{wInkw})):
\begin{equation}
{w_{(1)}}_{n;k}^{I}w_{(2)} = {w_{(1)}}_{n;k}^{I,0}w_{(2)}.
\end{equation}
}
\end{rema}

\begin{propo}\label{4.19}
Suppose that the $P(z)$-tensor product $(W_1\boxtimes_{P(z)} W_2,
Y_{P(z)}; \boxtimes_{P(z)})$ of $W_1$ and $W_2$ in ${\cal C}$ exists.
Then for any complex number $z_1\neq 0$, the $P(z_1)$-tensor product
of $W_1$ and $W_2$ in ${\cal C}$ also exists, and is given by
$(W_1\boxtimes_{P(z)} W_2, Y_{P(z)}; \boxtimes_{P(z_1)})$, where the
$P(z_1)$-intertwining map $\boxtimes_{P(z_1)}$ is defined by
\begin{equation}\label{tpzz_1}
\boxtimes_{P(z_1)}(w_{(1)}\otimes w_{(2)})=\sum_{n\in{\mathbb
C}}\sum_{k\in {\mathbb N}} {w_{(1)}}_{n;k}^{\boxtimes_{P(z)}}w_{(2)}
e^{\log z_1 (-n-1)}(\log z_1)^k
\end{equation}
for $w_{(1)}\in W_1$ and $w_{(2)}\in W_2$.
\end{propo}
\pf By Remark \ref{II1}, (\ref{tpzz_1}) indeed defines a
$P(z_1)$-product.  Given any $P(z_1)$-product $(W_3, Y_3; I_1)$ of
$W_1$ and $W_2$, let $I$ be the $P(z)$-product related to $I_1$ by
formula (\ref{log:zz_1}) with $I_1$, $I$ and $z_1$ in (\ref{log:zz_1})
replaced by $I$, $I_1$ and $z$, respectively, and with $p=0$.  Then {}from the
definition of $P(z)$-tensor product, there is a unique morphism $\eta$
{}from $(W_1\boxtimes_{P(z)} W_2, Y_{P(z)}; \boxtimes_{P(z)})$ to $(W_3,
Y_3; I)$.  Thus by (\ref{etaw1w2}) and (\ref{tpzz_1}) we
see that $\eta$ is also a morphism {}from the $P(z_1)$-product
$(W_1\boxtimes_{P(z)} W_2, Y_{P(z)}; \boxtimes_{P(z_1)})$ to $(W_3,
Y_3; I_1)$. The uniqueness of such a morphism follows similarly {}from
the uniqueness of a morphism {}from $(W_1\boxtimes_{P(z)} W_2, Y_{P(z)};
\boxtimes_{P(z)})$ to $(W_3, Y_3; I)$.  Hence $(W_1\boxtimes_{P(z)}
W_2, Y_{P(z)}; \boxtimes_{P(z_1)})$ is the $P(z_1)$-tensor product of
$W_1$ and $W_2$.  \epfv

\begin{rema}\label{intwmapdependsongeomdata} {\rm
In general, it will turn out that the existence of tensor product, and
the tensor product (generalized) module itself, do not depend on the
geometric data.  It is the intertwining map {}from the two modules to
the completion of their tensor product that encodes the geometric
information.}
\end{rema}

Generalizing Lemma 4.9 of \cite{tensor4}, we have:

\begin{propo}\label{span}
The generalized module $W_1\boxtimes_{P(z)}W_2$ (if it exists) is
spanned (as a vector space) by the (generalized-) weight components of
the elements of $\overline{W_1\boxtimes_{P(z)}W_2}$ of the form
$w_{(1)}\boxtimes_{P(z)} w_{(2)}$, for all $w_{(1)}\in W_1$ and
$w_{(2)}\in W_2$.
\end{propo}
\pf Denote by $W_0$ the vector subspace of $W_1\boxtimes_{P(z)}W_2$
spanned by all the weight components of all the elements of
$\overline{W_1\boxtimes_{P(z)} W_2}$ of the form
$w_{(1)}\boxtimes_{P(z)}w_{(2)}$ for $w_{(1)}\in W_1$ and $w_{(2)}\in
W_2$.  For a homogeneous vector $v\in V$ and arbitrary elements
$w_{(1)}\in W_1$ and $w_{(2)}\in W_2$, equating the
$x_0^{-1}x_1^{-m-1}$ coefficients of the Jacobi identity
(\ref{im:def}) gives
\begin{equation}\label{elm}
v_m({w_{(1)}}\boxtimes_{P(z)} w_{(2)})={w_{(1)}}\boxtimes_{P(z)} (v_m
w_{(2)})+\sum_{i\in {\mathbb N}}{m\choose i}z^{m-i}(v_i
w_{(1)})\boxtimes_{P(z)} w_{(2)}
\end{equation}
for all $m\in {\mathbb Z}$.  Note that the summation in the right-hand
side of (\ref{elm}) is always finite. Hence by taking arbitrary weight
components of (\ref{elm}) we see that $W_0$ is closed under the action
of $V$.  In case $V$ is M\"obius, a similar argument, using
(\ref{im:Lj}), shows that $W_0$ is stable under the action of
${\mathfrak s}{\mathfrak l}(2)$.  It is clear that $W_0$ is ${\mathbb
C}$-graded and $\tilde{A}$-graded.  Thus $W_0$ is a submodule of
$W_1\boxtimes_{P(z)} W_2$.

Now consider the quotient module
\[
W=(W_1\boxtimes_{P(z)} W_2)/W_0
\]
and let $\pi_{W}$ be the canonical map {}from $W_1\boxtimes_{P(z)}
W_2$ to $W$.  By the definition of $W_0$, we have
\[
\overline{\pi}_{W}\circ \boxtimes_{P(z)}=0,
\]
using the notation (\ref{actionofboxtensormap}).  The universal
property of the $P(z)$-tensor product then demands that $\pi_W=0$,
i.e., that $W_0=W_1\boxtimes_{P(z)} W_2$.

(Another argument: The image of the $P(z)$-intertwining map
$\boxtimes_{P(z)}$ lies in 
\[
\overline{W}_0 \subset \overline{W_1\boxtimes_{P(z)}W_2},
\]
so that $W_0$ is naturally a $P(z)$-product of $W_1$ and $W_2$, giving
rise to a (unique) $V$-module map
\[
f:W_1\boxtimes_{P(z)} W_2 \rightarrow W_0
\]
such that $\overline{f}$ takes each $w_{(1)}\boxtimes_{P(z)} w_{(2)}$
to $w_{(1)}\boxtimes_{P(z)} w_{(2)}$, by the universal property.
Writing
\[
\iota:W_0 \rightarrow W_1\boxtimes_{P(z)} W_2
\]
for the natural injection, we have that $\iota \circ f$ is the
identity map on $W_1\boxtimes_{P(z)} W_2$, by the universal property.
Thus $\iota$ is surjective (or, alternatively, $f$ is injective and is
$1$ on $W_0$), so that $W_0=W_1\boxtimes_{P(z)} W_2$.)
\epfv

It is clear {}from Definition \ref{pz-tp} that the tensor product
operation distributes over direct sums in the following sense:

\begin{propo}\label{tensorproductdistributes}
For $U_1, \dots, U_{k}$, $W_1, \dots, W_{l}\in \ob{\cal C}$,
suppose that each $U_{i}\boxtimes_{P(z)}W_{j}$ exists. Then
$(\coprod_{i}U_{i})\boxtimes_{P(z)}(\coprod_{j}W_{j})$ exists and
there is a natural isomorphism
\[
\biggl(\coprod_{i}U_{i}\biggr)\boxtimes_{P(z)}
\biggl(\coprod_{j}W_{j}\biggr)\stackrel{\sim}
{\rightarrow} \coprod_{i,j}U_{i}\boxtimes_{P(z)}W_{j}.\hspace{4em}\square
\]
\end{propo}

\begin{rema}\label{bifunctor} {\rm
It is of course natural to view the $P(z)$-tensor product as a
bifunctor: Suppose that ${\cal C}$ is a full subcategory of either
${\cal M}_{sg}$ or ${\cal GM}_{sg}$ (recall Notation \ref{MGM}) such
that for all $W_{1}, W_{2}\in \ob{\cal C}$, the $P(z)$-tensor product
of $W_{1}$ and $W_{2}$ exists in ${\cal C}$.  Then $\boxtimes_{P(z)}$
provides a (bi)functor
\begin{equation}
\boxtimes_{P(z)}:{\cal C} \times  {\cal C} \rightarrow {\cal C}
\end{equation}
as follows: For $W_{1}, W_{2}\in \ob{\cal C}$,
\begin{equation}
\boxtimes_{P(z)}(W_{1}, W_{2}) = W_1\boxtimes_{P(z)} W_2 \in \ob{\cal C}
\end{equation}
and for $V$-module maps
\begin{equation}
\sigma_1 : W_1 \rightarrow W_3,
\end{equation}
\begin{equation}
\sigma_2 : W_2 \rightarrow W_4
\end{equation}
with $W_{3}, W_{4}\in \ob{\cal C}$, we have the $V$-module map,
denoted
\begin{equation}
\boxtimes_{P(z)}(\sigma_{1}, \sigma_{2}) = \sigma_1\boxtimes_{P(z)}
\sigma_2,
\end{equation}
{}from $W_1\boxtimes_{P(z)} W_2$ to $W_3\boxtimes_{P(z)} W_4$, defined
by the universal property of the $P(z)$-tensor product
$W_1\boxtimes_{P(z)} W_2$ and the fact that the composition of
$\boxtimes_{P(z)}$ with $\sigma_1 \otimes \sigma_2$ is a
$P(z)$-intertwining map
\begin{equation}
\boxtimes_{P(z)} \circ (\sigma_{1} \otimes \sigma_{2}):W_1\otimes W_2
\rightarrow \overline{W_3\boxtimes_{P(z)} W_4}.
\end{equation}
Note that it is the effect of this bifunctor on morphisms (rather than
on objects) that exhibits the role of the geometric data.
}
\end{rema}

We obtain right exact functors by fixing one of the generalized
modules in Remark \ref{bifunctor}\footnote{We thank Ingo Runkel for
asking us whether our tensor product functors are right exact.}:

\begin{propo}
In the setting of Remark \ref{bifunctor}, for $W \in \ob{\cal C}$ the
functors $W\boxtimes_{P(z)} \cdot$ and $\cdot \boxtimes_{P(z)} W$ are
right exact.
\end{propo}
\pf Let
\[
W_1 \stackrel{\sigma_1}{\longrightarrow} W_2
\stackrel{\sigma_2}{\longrightarrow} W_3 \longrightarrow 0
\]
be exact in ${\cal C}$.  We show that
\[
W \boxtimes_{P(z)} W_1 
\stackrel{1_W \boxtimes_{P(z)} \sigma_1}{\longrightarrow}
W \boxtimes_{P(z)} W_2
\stackrel{1_W \boxtimes_{P(z)} \sigma_2}{\longrightarrow}
W \boxtimes_{P(z)} W_3
\longrightarrow 0
\]
is exact; the proof of right exactness for $\cdot \boxtimes_{P(z)} W$
is completely analogous.

For the surjectivity of ${1_W \boxtimes \sigma_2}$, we observe that
the elements
\[
(1_W \boxtimes \sigma_2)(\pi_n (w \boxtimes w_{(2)}))
\]
for $w \in W$, $w_{(2)} \in W_2$ and $n \in {\mathbb C}$ span $W
\boxtimes W_3$ (we are dropping the subscripts $P(z)$), since this
element equals
\[
\pi_n (\overline{1_W \boxtimes \sigma_2}(w \boxtimes w_{(2)}))
=\pi_n (w \boxtimes \sigma_2 (w_{(2)})),
\]
and these elements span $W \boxtimes W_3$ by the surjectivity of
$\sigma_2$ and Proposition \ref{span}.

Since 
\[
(1_W \boxtimes \sigma_2)(1_W \boxtimes \sigma_1)
= 1_W \boxtimes \sigma_2\sigma_1 = 0,
\]
it remains only to show that the natural (surjective) module map
\[
\theta:(W_1 \boxtimes W_2)/{\rm{Im}}\,(1_W \boxtimes \sigma_1)
\rightarrow W_1 \boxtimes W_3
\]
is injective.  Noting that
\[
\overline{(W_1 \boxtimes W_2)/{\rm{Im}}\,(1_W \boxtimes \sigma_1)}
= \overline{(W_1 \boxtimes W_2)}/
\overline{{\rm{Im}}\,(1_W \boxtimes \sigma_1)},
\]
we characterize $\theta$ by:
\[
\overline{\theta}(w \boxtimes w_{(2)}
+ \overline{{\rm{Im}}\,(1_W \boxtimes \sigma_1)})
= \overline{1_W \boxtimes \sigma_2}(w \boxtimes w_{(2)}).
\]
We construct a $P(z)$-intertwining map
\[
I:W \otimes W_3 \rightarrow
\overline{(W \boxtimes W_2)/{\rm{Im}}\,(1_W \boxtimes \sigma_1)}
\]
as follows:  For $w \in W$ and $w_{(3)} \in W_3$ set
\[
I(w \otimes w_{(3)}) = w \boxtimes w_{(2)}
+ \overline{{\rm{Im}}\,(1_W \boxtimes \sigma_1)}
\]
where $w_{(2)} \in W_2$ is such that
\[
\sigma_2 (w_{(2)})=w_{(3)}.
\]
Then $I$ is well defined because for $w'_{(2)} \in W_2$ with
$\sigma_2 (w'_{(2)})=w_{(3)}$,
\[
w \otimes (w_{(2)}-w'_{(2)}) \in w \boxtimes {\rm{Ker}}\, \sigma_2
= w \boxtimes {\rm{Im}}\, \sigma_1 \subset
\overline{{\rm{Im}}\, (1_W \otimes \sigma_1)},
\]
and it is straightforward to verify that $I$ is in fact a
$P(z)$-intertwining map.  Thus we have a module map
\[
\eta: W \boxtimes W_3 \rightarrow (W \boxtimes W_2)/ 
{\rm{Im}}\, (1_W \otimes \sigma_1)
\]
such that
\[
\overline{\eta}(w \boxtimes w_{(3)})
= w \boxtimes w_{(2)} + \overline{{\rm{Im}}\, (1_W \otimes \sigma_1)},
\]
with the elements as above.  Then
\begin{eqnarray*}
\lefteqn{\overline{\eta \circ \theta}(w \boxtimes w_{(2)}
+ \overline{{\rm{Im}}\, (1_W \otimes \sigma_1)})
= {\overline{\eta}}(\overline{1_W \otimes \sigma_2}(w \boxtimes
w_{(2)}))}\\
&& \quad\quad\quad\quad\quad\quad\quad\quad\quad\quad
= {\overline{\eta}}(w \boxtimes \sigma_2 (w_{(2)}))\\
&& \quad\quad\quad\quad\quad\quad\quad\quad\quad\quad
= {\overline{\eta}}(w \boxtimes w_{(3)})\\
&& \quad\quad\quad\quad\quad\quad\quad\quad\quad\quad
= w \boxtimes w_{(2)} + \overline{{\rm{Im}}\, (1_W \otimes \sigma_1)},
\end{eqnarray*}
which shows that $\eta \circ \theta$ is the identity map, and so
$\theta$ is injective, as desired.
\epfv

We now discuss the simplest examples of $P(z)$-tensor products---those
in which one or both of $W_{1}$ or
$W_{2}$ is $V$ itself (viewed as a (generalized) $V$-module); 
we suppose here that $V\in \ob \mathcal{C}$. Since the discussion of the 
case in which both
$W_{1}$ and $W_{2}$ are $V$ turns out to be no simpler than the case 
in which $W_{1}=V$,
we shall discuss only the two more general 
cases $W_{1}=V$ and $W_{2}=V$. 

\begin{exam}\label{expl-vw}
{\rm Let $(W, Y_{W})$ be an object of $\mathcal{C}$. 
The vertex operator map $Y_{W}$ gives a $P(z)$-intertwining 
map 
\[
I_{Y_{W}, p}=Y_{W}(\cdot, z)\cdot: V\otimes W\to \overline{W}
\]
for any fixed $p\in \Z$ (recall Proposition \ref{im:correspond} and Remark 
\ref{mod-sub}). We claim that 
$(W, Y_{W}; Y_{W}(\cdot, z)\cdot)$ is the $P(z)$-tensor product of $V$ and 
$W$ in $\mathcal{C}$. In fact, let $(W_{3}, Y_{3}; I)$ be a $P(z)$-product 
of $V$ and $W$ in $\mathcal{C}$ and suppose that there exists 
a module map $\eta: W\to W_{3}$ such that 
\begin{equation}\label{v-tensor-w-1}
\overline{\eta}\circ (Y_{W}(\cdot, z)\cdot)=I.
\end{equation}
Then for $w\in W$, we must have 
\begin{eqnarray}\label{v-tensor-w-2}
\eta(w)&=&\eta(Y_{W}(\mathbf{1}, z)w)\nn
&=&(\overline{\eta}\circ (Y_{W}(\cdot, z)\cdot))(\mathbf{1}\otimes w)\nn
&=&I(\mathbf{1}\otimes w),
\end{eqnarray}
so that $\eta$ is unique if it exists.
We now define $\eta: W\to \overline{W}_{3}$ using (\ref{v-tensor-w-2}). 
We shall show that $\eta(W)\subset W_{3}$ and that $\eta$ has the 
desired properties.
Since 
$I$ is a $P(z)$-intertwining map of type ${W_{3}\choose VW}$, it corresponds 
to a logarithmic intertwining 
operator $\mathcal{Y}=\mathcal{Y}_{I, p}$ of the same type, according to 
Proposition \ref{im:correspond}. Since $L(-1)\mathbf{1}=0$, 
we have 
\[
\frac{d}{dx}\mathcal{Y}(\mathbf{1}, x)=\mathcal{Y}(L(-1)\mathbf{1}, x)=0.
\]
Thus $\mathcal{Y}(\mathbf{1}, x)$ 
is simply the constant map $\mathbf{1}_{-1; 0}^{\mathcal{Y}}: W\to W_{3}$ 
(using the notation (\ref{log:map})), and this map preserves (generalized) 
weights, by Proposition \ref{log:logwt}(b).  By
Proposition 
\ref{im:correspond}, $I=I_{\mathcal{Y}, p}$, so that
\begin{eqnarray*}
\eta(w)&=&I(\mathbf{1}\otimes w)\nn
&=&I_{\mathcal{Y}, p}(\mathbf{1}\otimes w)\nn
&=&\mathbf{1}_{-1; 0}^{\mathcal{Y}}w
\end{eqnarray*}
for $w\in W$.
So $\eta=\mathbf{1}_{-1; 0}^{\mathcal{Y}}$ is a linear map {}from 
$W$ to $W_{3}$ preserving (generalized) weights. Using the 
Jacobi identity (\ref{im:def}) for the $P(z)$-intertwining map $I$ and the 
fact that $Y(u, x_{0})\mathbf{1}\in V[[x_{0}]]$ for $u\in V$, we obtain
\begin{eqnarray*}
\eta(Y_{W}(u, x)w)&=&I(\mathbf{1} \otimes Y_{W}(u, x)w)\nn
&=&Y_{3}(u, x)I(\mathbf{1} \otimes w)
-\res_{x_{0}}z^{-1}\delta\left(\frac{x-x_{0}}{z}\right)
I(Y(u, x_{0})\mathbf{1} \otimes w)\nn
&=&Y_{3}(u, x)I(\mathbf{1} \otimes w)\nn
&=&Y_{3}(u, x)\eta(w)
\end{eqnarray*}
for $u\in V$ and $w\in W$, proving that
$\eta$ is a module map when $V$ is a conformal vertex algebra, and
when $V$ is M\"obius, $\eta$ also commutes with the action of
${\mathfrak s}{\mathfrak l}(2)$, by (\ref{im:Lj}).  For $w\in W$, 
\begin{eqnarray}\label{v-tensor-w-3}
(\overline{\eta}\circ (Y_{W}(\cdot, z)\cdot))(\mathbf{1}\otimes w)&=&
\overline{\eta}(Y_{W}(\mathbf{1}, z)w)\nn
&=&\eta(w)\nn
&=&I(\mathbf{1}\otimes w).
\end{eqnarray}
Using the Jacobi identity for $P(z)$-intertwining maps, we obtain
\begin{eqnarray}\label{int-recurrence-rel}
\lefteqn{I(Y(u, x_{0})v\otimes w)}\nn
&&=\res_{x}x_{0}^{-1}\delta\left(\frac{x-z}{x_{0}}\right)
Y_3(u, x)I(v\otimes w)-\res_{x}x_{0}^{-1}\delta\left(\frac{z-x}{-x_{0}}\right)
I(v\otimes Y_W(u, x)w) \;\;\;\;\;\;\;
\end{eqnarray}
for $u, v\in V$ and $w\in W$. 
Since $\eta$ is a module map and $Y_{W}(\cdot, z)\cdot$ is a $P(z)$-intertwining 
map of type ${W\choose VW}$, $\overline{\eta}\circ Y_{W}(\cdot, z)\cdot$
is a $P(z)$-intertwining map of  type ${W_{3}\choose VW}$. In particular, 
(\ref{int-recurrence-rel}) holds when we replace $I$ by 
$\overline{\eta}\circ Y_{W}(\cdot, z)\cdot$. Using
(\ref{int-recurrence-rel}) for $v=\mathbf{1}$ together with (\ref{v-tensor-w-3}), 
we obtain
\[
(\overline{\eta}\circ (Y_{W}(\cdot, z)\cdot))(u\otimes w)=I(u\otimes w)
\]
for $u\in V$ and $w\in W$, proving (\ref{v-tensor-w-1}), as desired.
Thus $(W, Y_{W}; Y_{W}(\cdot, z)\cdot)$ is the $P(z)$-tensor product
of $V$ and $W$ in $\mathcal{C}$. }
\end{exam}

\begin{exam}\label{expl-wv}
{\rm Let $(W, Y_{W})$ be an object of $\mathcal{C}$. In order to
construct the $P(z)$-tensor product $W\boxtimes_{P(z)} V$, recall {}from
(\ref{Omega_r}) and Proposition \ref{log:omega} that
$\Omega_{p}(Y_{W})$ is a logarithmic intertwining operator of type
${W\choose WV}$. It involves only integral powers of the formal
variable and no logarithms, and it is independent of $p$.  In fact,
\[
\Omega_{p}(Y_{W})(w, x)v=e^{xL(-1)}Y_{W}(v, -x)w
\]
for $v\in V$ and $w\in W$. For $q\in \Z$, 
\[
I_{\Omega_{p}(Y_{W}), q}
=\Omega_{p}(Y_{W})(\cdot, z)\cdot: W\otimes 
V\to \overline{W}
\]
is a $P(z)$-intertwining map of the same type and is independent of
$q$.  We claim that $(W, Y_{W}; \Omega_{p}(Y_{W})(\cdot, z)\cdot)$ is
the $P(z)$-tensor product of $W$ and $V$ in $\mathcal{C}$. In fact,
let $(W_{3}, Y_{3}; I)$ be a $P(z)$-product of $W$ and $V$ in
$\mathcal{C}$ and suppose that there exists a module map $\eta: W\to
W_{3}$ such that
\begin{equation}\label{w-tensor-v-1}
\overline{\eta}\circ \Omega_{p}(Y_{W})(\cdot, z)\cdot=I.
\end{equation}
For $w\in W$, we must have 
\begin{eqnarray}\label{w-tensor-v-3}
\eta(w)&=&\eta(Y_{W}(\mathbf{1}, -z)w)\nn
&=&e^{-zL(-1)}\overline{\eta}(e^{zL(-1)}Y_{W}(\mathbf{1}, -z)w)\nn
&=&e^{-zL(-1)}\overline{\eta}(\Omega_{p}(Y_{W})(w, z)\mathbf{1}))\nn
&=&e^{-zL(-1)}(\overline{\eta}\circ (\Omega_{p}(Y_{W})(\cdot, z)\cdot))
(w \otimes \mathbf{1})\nn
&=&e^{-zL(-1)}I(w\otimes \mathbf{1}),
\end{eqnarray}
and so $\eta$ is unique if it exists.  (Note that the right-hand side
of (\ref{w-tensor-v-3}) is indeed defined, in view of (\ref{im:ltc}).)
We now define $\eta: W\to \overline{W}_{3}$ by (\ref{w-tensor-v-3}).
Consider the logarithmic intertwining operator $\mathcal{Y}=
\mathcal{Y}_{I, q}$ that corresponds to $I$ by Proposition  
\ref{im:correspond}. 
Using Proposition \ref{im:correspond},
(\ref{branch1})--(\ref{log:fsub}),
(\ref{log:subs})
and the equality
\begin{eqnarray*}
l_{q}(-z)&=&\log |-z|+i(\arg (-z)+2\pi q)\nn
&=&\left\{\begin{array}{ll}
\log |z|+i(\arg z+\pi+2\pi q),&0\le \arg z< \pi\\
\log |z|+i(\arg z-\pi+2\pi q),&\pi\le \arg z< 2\pi
\end{array}\right.\nn
&=&\left\{\begin{array}{ll}
l_{q}(z)+\pi i,&0\le \arg z< \pi\\
l_{q}(z)-\pi i,&\pi\le \arg z< 2\pi,
\end{array}\right.\nn
\end{eqnarray*}
we have
\begin{eqnarray*}
e^{-zL(-1)}I(w\otimes \mathbf{1})
&=&e^{-zL(-1)}\mathcal{Y}(w, e^{l_{q}(z)})\mathbf{1}\nn
&=&e^{-xL(-1)}\mathcal{Y}(w, x)\mathbf{1}|_{x^{n}=e^{nl_{q}(z)},\; 
(\log x)^{m}=(l_{q}(z))^{m},\; n\in \C, \;m\in \N}\nn
&=&e^{yL(-1)}\mathcal{Y}(w, e^{\pm \pi i}y)\mathbf{1}
|_{y^{n}=e^{nl_{q}(-z)},\; (\log y)^{m}=(l_{q}(-z))^{m},\;
n\in \C, \; m\in \N},
\end{eqnarray*}
where $e^{\pm \pi i}$ is $e^{-\pi i}$ when $0\le \arg z<\pi$ and is 
$e^{\pi i}$ when $\pi\le \arg z<2 \pi$. Then by (\ref{Omega_r}),
we see that 
$\eta(w)=e^{-zL(-1)}I(w\otimes \mathbf{1})$ is equal to 
$\Omega_{-1}(\mathcal{Y})(\mathbf{1}, e^{l_{q}(-z)})w$ 
when $0\le \arg z<\pi$ and 
is equal to $\Omega_{0}(\mathcal{Y})(\mathbf{1}, e^{l_{q}(-z)})w$ when 
$\pi\le \arg z<2 \pi$. By Proposition \ref{log:omega}, 
$\Omega_{-1}(\mathcal{Y})$ and $\Omega_{0}(\mathcal{Y})$ are
logarithmic intertwining operators of type ${W_{3}\choose VW}$. 
As in Example \ref{expl-vw}, we see that 
$\Omega_{-1}(\mathcal{Y})(\mathbf{1}, y)$ 
and $\Omega_{0}(\mathcal{Y})(\mathbf{1}, y)$ are equal to 
$\mathbf{1}_{-1, 0}^{\Omega_{-1}(\mathcal{Y})}$ and
$\mathbf{1}_{-1, 0}^{\Omega_{0}(\mathcal{Y})}$, respectively, and 
these maps preserve (generalized) weights. Therefore $\eta$ 
is a linear map {}from $W$ to $W_{3}$ preserving (generalized) weights.
Using the Jacobi identity (\ref{im:def})
for the $P(z)$-intertwining map $I$ and the 
fact that $Y(u, x_{1})\mathbf{1}\in V[[x_{1}]]$, we have
\begin{eqnarray*}
\eta(Y_{W}(u, x_{0})w)&=&e^{-zL(-1)}I(Y_{W}(u, x_{0})w\otimes \mathbf{1})\nn
&=&\res_{x_{1}}x_{0}^{-1}\delta\left(\frac{x_{1}-z}{x_{0}}\right)
e^{-zL(-1)}Y_{3}(u, x_{1})I(w\otimes \mathbf{1})\nn
&& -\res_{x_{1}}x_{0}^{-1}\delta\left(\frac{z-x_{1}}{-x_{0}}\right)
e^{-zL(-1)}I(w\otimes Y(u, x_{1})\mathbf{1})\nn
&=&e^{-zL(-1)}Y_{3}(u, x_{0}+z)I(w\otimes \mathbf{1})\nn
&=&Y_{3}(u, x_{0})e^{-zL(-1)}I(w\otimes \mathbf{1})\nn
&=&Y_{3}(u, x_{0})\eta(w)
\end{eqnarray*}
for $u\in V$ and $w\in W$, proving that $\eta$ is a module map when
$V$ is a conformal vertex algebra.  As in Example \ref{expl-vw}, when
$V$ is M\"obius, $\eta$ also commutes with the action of 
${\mathfrak s}{\mathfrak l}(2)$, this time by (\ref{im:Lj2}) together with
(\ref{log:SL2-1}) with $x$ specialized to $-z$.  For $w\in W$, 
\begin{eqnarray}\label{w-tensor-v-4}
(\overline{\eta}\circ (\Omega_{p}(Y_{W})(\cdot, z)\cdot))(w \otimes \mathbf{1})
&=&\overline{\eta}(e^{zL(-1)}Y_{W}(\mathbf{1}, -z)w)\nn
&=&e^{zL(-1)}\eta(w)\nn
&=&e^{zL(-1)}e^{-zL(-1)}I(w\otimes \mathbf{1})\nn
&=&I(w\otimes \mathbf{1}).
\end{eqnarray}
Since both $\overline{\eta}\circ (\Omega_{p}(Y_{W})(\cdot, z)\cdot)$ and $I$ 
are $P(z)$-intertwining maps of type ${W_{3}\choose WV}$,
using the Jacobi identity for $P(z)$-intertwining operators
and (\ref{w-tensor-v-4}) (cf. Example \ref{expl-vw}), we have
\[
(\overline{\eta}\circ (\Omega_{p}(Y_{W})(\cdot, z)\cdot))
(w\otimes v)=I(w\otimes v)
\]
for $v\in V$ and $w\in W$, proving (\ref{w-tensor-v-1}). 
Thus $(W, Y_{W}; \Omega_{p}(Y_{W})(\cdot, z)\cdot)$  is the 
$P(z)$-tensor product of $W$ and 
$V$ in $\mathcal{C}$.}
\end{exam}

We discussed the important special class of finitely reductive vertex
operator algebras in the Introduction.  In case $V$ is a finitely
reductive vertex operator algebra, the $P(z)$-tensor product always
exists, as we are about to establish (following \cite{tensor1} and
\cite{tensor3}).  As in the Introduction, the definition of finite
reductivity is:

\begin{defi}\label{finitelyreductive}{\rm
A vertex operator algebra $V$ is {\it finitely reductive} if
\begin{enumerate}
\item Every $V$-module is completely reducible.
\item There are only finitely many irreducible $V$-modules (up to
equivalence).
\item All the fusion rules (the dimensions of the spaces of
intertwining operators among triples of modules) for $V$ are finite.
\end{enumerate}
}
\end{defi}

\begin{rema} {\rm
In this case, every $V$-module is of course a {\it finite} direct sum
of irreducible modules.  Also, the third condition holds if the
finiteness of the fusion rules among triples of only {\it irreducible}
modules is assumed.
}
\end{rema}

\begin{rema} {\rm
We are of course taking the notion of $V$-module so that the grading
restriction conditions are the ones described in Remark
\ref{moduleswiththetrivialgroup}, formulas (\ref{Wn+k=0}) and
(\ref{dimWnfinite}); in particular, $V$-modules are understood to be
$\C$-graded.  Recall {}from Remark \ref{congruent} that for an
irreducible module, all its weights are congruent to one another
modulo $\Z$.  Thus for an irreducible module, our grading-truncation
condition (\ref{Wn+k=0}) amounts exactly to the condition that the
real parts of the weights are bounded {}from below.  In
\cite{tensor1}--\cite{tensor3}, boundedness of the real parts of the
weights {}from below was our grading-truncation condition in the
definition of the notion of module for a vertex operator algebra.
Thus the first two conditions in the notion of finite reductivity are
the same whether we use the current grading restriction conditions in
the definition of the notion of module or the corresponding conditions
in \cite{tensor1}--\cite{tensor3}.  As for intertwining operators,
recall {}from Remark \ref{ordinaryandlogintwops} and Corollary
\ref{powerscongruentmodZ} that when the first two conditions are
satisfied, the notion of (ordinary, non-logarithmic) intertwining
operator here coincides with that in \cite{tensor1} because the
truncation conditions agree.  Also, in this setting, by Remark
\ref{log:ordi}, the logarithmic and ordinary intertwining operators
are the same, and so the spaces of intertwining operators ${\cal
V}^{W_3}_{W_1\,W_2}$ and fusion rules $N^{W_3}_{W_1\,W_2}$ in
Definition \ref{fusionrule} have the same meanings as in
\cite{tensor1}.  Thus the notion of finite reductivity for a vertex
operator algebra is the same whether we use the current grading
restriction and truncation conditions in the definitions of the
notions of module and of intertwining operator or the corresponding
conditions in \cite{tensor1}--\cite{tensor3}.  In particular, finite
reductivity of $V$ according to Definition \ref{finitelyreductive} is
equivalent to the corresponding notion, ``rationality'' (recall the
Introduction) in \cite{tensor1}--\cite{tensor3}.  }
\end{rema}

\begin{rema} {\rm
For a vertex operator algebra $V$ (in particular, a finitely reductive
one), the category ${\cal M}$ of $V$-modules coincides with the
category ${\cal M}_{sg}$ of strongly graded $V$-modules; recall
Notation \ref{MGM}.
}
\end{rema}

For the rest of Section 4.1, let us assume that $V$ is a finitely
reductive vertex operator algebra.  We shall now show that
$P(z)$-tensor products always exist in the category ${\cal M}={\cal
M}_{sg}$ of $V$-modules, in the sense of Definition \ref{pz-tp}.

The considerations from here through (\ref{fusionrulerelation}) also
hold, with natural adjustments, for finite-dimensional modules for a
semisimple Lie algebra (even though there are infinitely many
irreducible modules up to equivalence) or for a finite group or for a
compact group, etc., but in such classical contexts, one does not
ordinarily express things in this way because one knows {\it a priori}
that the tensor product functors exist and satisfy natural
associativity as in (\ref{calWassociativity}), (\ref{wassociativity}).
What we do now shows how to build tensor product functors with
knowledge ``only'' of the spaces of intertwining maps, and uses this
to motivate how to approach the problem of constructing appropriate
natural associativity isomorphisms, whether or not our vertex algebra
$V$ is a finitely reductive vertex operator algebra.

Consider $V$-modules $W_{1}$, $W_{2}$ and $W_{3}$.  We know that
\begin{eqnarray}
N^{W_3}_{W_1\,W_2} = \dim {\cal V}^{W_3}_{W_1\,W_2} < \infty
\end{eqnarray}
and {}from Proposition \ref{im:correspond}, we also have
\begin{eqnarray}
N^{W_3}_{W_1\,W_2} = \dim {\cal M}[P(z)]_{W_{1}W_{2}}^{W_{3}} = \dim
{\cal M}_{W_{1}W_{2}}^{W_{3}} < \infty
\end{eqnarray}
(recall Definition \ref{im:imdef}).

The natural evaluation map
\begin{eqnarray}
W_{1}\otimes W_{2}\otimes {\cal M}^{W_{3}}_{W_{1}W_{2}}&\to& \overline{W}_{3}
\nno\\
w_{(1)}\otimes w_{(2)}\otimes I&\mapsto& I(w_{(1)}\otimes w_{(2)})
\end{eqnarray}
gives a natural map
\begin{equation}
{\cal F}[P(z)]^{W_{3}}_{W_{1}W_{2}}: W_{1}\otimes W_{2}\to 
\mbox{\rm Hom}({\cal M}_{W_{1}W_{2}}^{W_{3}}, \overline{W}_{3})=
({\cal M}^{W_{3}}_{W_{1}W_{2}})^{*}\otimes \overline{W}_{3}.
\end{equation}
Since $\dim {\cal M}_{W_{1}W_{2}}^{W_{3}} < \infty$, $({\cal
M}^{W_{3}}_{W_{1}W_{2}})^{*}\otimes W_{3}$ is a $V$-module (with
finite-dimensional weight spaces) in the obvious way, and the map
${\cal F}[P(z)]^{W_{3}}_{W_{1}W_{2}}$ is clearly a $P(z)$-intertwining
map, where we make the identification
\begin{equation}
({\cal M}^{W_{3}}_{W_{1}W_{2}})^{*}\otimes \overline{W}_{3}
=\overline{({\cal M}^{W_{3}}_{W_{1}W_{2}})^{*}\otimes W_{3}}.
\end{equation}
This gives us a natural $P(z)$-product for the category ${\cal M} =
{\cal M}_{sg}$ (recall Definition \ref{pz-product}).  Moreover, we
have a natural linear injection
\begin{eqnarray}
i: {\cal M}^{W_{3}}_{W_{1}W_{2}}&\to &
\mbox{\rm Hom}_{V}(({\cal M}^{W_{3}}_{W_{1}W_{2}})^{*}\otimes W_{3}, W_{3})\nno\\
I&\mapsto &(f\otimes w_{(3)}\mapsto f(I)w_{(3)})
\end{eqnarray}
which is an isomorphism if $W_{3}$ is irreducible, since in this
case, 
\[
\mbox{\rm Hom}_{V}(W_{3}, W_{3})\simeq {\C}
\]
(see \cite{FHL}, Remark 4.7.1).  On the other hand, the natural map
\begin{eqnarray}
h:\mbox{\rm Hom}_{V}(({\cal M}^{W_{3}}_{W_{1}W_{2}})^{*}
\otimes W_{3}, W_{3})&\to &
{\cal M}^{W_{3}}_{W_{1}W_{2}}\nno\\
\eta&\mapsto &\overline{\eta}\circ {\cal F}[P(z)]^{W_{3}}_{W_{1}W_{2}}
\end{eqnarray}
given by composition clearly satisfies the condition that
\begin{equation}\label{hiI=I}
h(i(I))=I,
\end{equation}
so that if $W_{3}$ is irreducible, the maps $h$ and $i$ are mutually inverse 
isomorphisms and we have the property that for any $I\in 
{\cal M}^{W_{3}}_{W_{1}W_{2}}$, there exists a unique $\eta$ such that 
\begin{equation}\label{I=etabarF}
I=\overline{\eta}\circ {\cal F}[P(z)]^{W_{3}}_{W_{1}W_{2}}
\end{equation}
(cf. Definition \ref{pz-tp}).

Using this, we can now show, in the next result, that $P(z)$-tensor
products always exist for the category of modules for a finitely
reductive vertex operator algebra, and we shall in fact exhibit the
$P(z)$-tensor product.  Note that there is no need to assume that
$W_{1}$ and $W_{2}$ are irreducible in the formulation or proof, but
by Proposition \ref{tensorproductdistributes}, the case in which
$W_{1}$ and $W_{2}$ are irreducible is in fact sufficient, and the
tensor product operation is canonically described using only the
spaces of intertwining maps among triples of {\it irreducible}
modules.

\begin{propo}\label{construcofPztensorprod-finredcase}
Let $V$ be a finitely reductive vertex operator algebra and let
$W_{1}$ and $W_{2}$ be $V$-modules.  Then
$(W_{1}\boxtimes_{P(z)}W_{2}, Y_{P(z)}; \boxtimes_{P(z)})$ exists, and
in fact
\begin{equation}\label{Pztensorprodfinitelyredcase}
W_{1}\boxtimes_{P(z)}W_{2}=\coprod_{i=1}^{k}
({\cal M}^{M_{i}}_{W_{1}W_{2}})^{*}\otimes M_{i},
\end{equation}
where $\{ M_{1}, \dots, M_{k}\}$ is a set of representatives of the
equivalence classes of irreducible $V$-modules, and the right-hand
side of (\ref{Pztensorprodfinitelyredcase}) is equipped with the
$V$-module and $P(z)$-product structure indicated above. That is,
\begin{equation}
\boxtimes_{P(z)}=\sum_{i=1}^{k}{\cal F}[P(z)]^{M_{i}}_{W_{1}W_{2}}.
\end{equation}
\end{propo}
\pf
{From} the comments above and the definitions, it is clear that we have a 
$P(z)$-product. Let $(W_{3}, Y_{3}; I)$ be any $P(z)$-product. Then 
$W_{3}=\coprod_{j}U_{j}$ where $j$ ranges through a finite set and each
$U_{j}$ is irreducible. Let $\pi_{j}: W_{3}\to U_{j}$ denote the $j$-th 
projection. A module map $\eta:\coprod_{i=1}^{k}
({\cal M}^{M_{i}}_{W_{1}W_{2}})^{*}\otimes M_{i}\to W_{3}$ amounts to 
module maps
$$\eta_{ij}:  ({\cal M}^{M_{i}}_{W_{1}W_{2}})^{*}\otimes M_{i}\to U_{j}$$
for each $i$ and $j$ such that $U_{j}\simeq M_{i}$, and 
$I=\overline{\eta}\circ \boxtimes_{P(z)}$ if and only if
\[
\overline{\pi}_{j}\circ I=\overline{\eta}_{ij}\circ 
{\cal F}^{M_{i}}_{W_{1}W_{2}}
\]
for each $i$ and $j$, the bars having the obvious meaning. But 
$\overline{\pi}_{j}\circ I$ is a $P(z)$-intertwining map of type 
${U_{j}}\choose {W_{1}W_{2}}$, and so 
$\overline{\iota}\circ \overline{\pi}_{j}\circ I\in 
{\cal M}^{M_{i}}_{W_{1}W_{2}}$,
where $\iota: U_{j}\stackrel{\sim}{\to}M_{i}$ is a fixed isomorphism. 
Denote this map by $\tau$. Thus what we finally want is a unique module map
\[
\theta: ({\cal M}^{M_{i}}_{W_{1}W_{2}})^{*}\otimes M_{i}\to M_{i}
\]
such that 
\[
\tau=\overline{\theta}\circ {\cal F}[P(z)]^{M_{i}}_{W_{1}W_{2}}.
\]
But we in fact have such a unique $\theta$, by
(\ref{hiI=I})--(\ref{I=etabarF}).  \epfv

\begin{rema} {\rm
By combining Proposition \ref{construcofPztensorprod-finredcase} with
Proposition \ref{im:correspond}, we can express $W_{1}\boxtimes_{P(z)}
W_{2}$ in terms of ${\cal V}^{M_{i}}_{W_{1}W_{2}}$ in place of ${\cal
M}^{M_{i}}_{W_{1}W_{2}}$.
}
\end{rema}

\begin{rema} {\rm
If we know the fusion rules among triples of irreducible $V$-modules,
then {}from Proposition \ref{construcofPztensorprod-finredcase} we know
all the $P(z)$-tensor product modules, up to equivalence; that is, we
know the multiplicity of each irreducible $V$-module in each
$P(z)$-tensor product module.  But recall that the $P(z)$-tensor
product structure of $W_{1}\boxtimes_{P(z)} W_{2}$ involves much more
than just the $V$-module structure.
}
\end{rema}

As we discussed in the Introduction, the main theme of this work is to
construct natural ``associativity'' isomorphisms between triple tensor
products of the shape $W_{1}\boxtimes (W_{2}\boxtimes W_{3})$ and
$(W_{1}\boxtimes W_{2})\boxtimes W_{3}$, for (generalized) modules
$W_{1}$, $W_{2}$ and $W_{3}$.  In the finitely reductive case, let
$W_{1}$, $W_{2}$ and $W_{3}$ be $V$-modules.  By Proposition
\ref{construcofPztensorprod-finredcase}, we have, as $V$-modules,
\begin{eqnarray}\label{W1(W2W3)}
\lefteqn{W_{1}\boxtimes_{P(z)}(W_{2}\boxtimes_{P(z)}W_{3})
= W_{1}\boxtimes_{P(z)}\left(\coprod_{i=1}^{k} M_{i} \otimes ({\cal
M}^{M_{i}}_{W_{2}W_{3}})^{*}\right)}
\nno\\
&&
=\coprod_{i=1}^{k} (W_{1}\boxtimes_{P(z)} M_{i}) \otimes ({\cal
M}^{M_{i}}_{W_{2}W_{3}})^{*}
\nno\\
&&
=\coprod_{i=1}^{k} \left(\coprod_{j=1}^{k} ({\cal
M}^{M_{j}}_{W_{1}M_{i}})^{*} \otimes M_{j}\right) \otimes ({\cal
M}^{M_{i}}_{W_{2}W_{3}})^{*}
\nno\\
&&
=\coprod_{j=1}^{k} \left(\coprod_{i=1}^{k} ({\cal
M}^{M_{j}}_{W_{1}M_{i}})^{*} \otimes ({\cal
M}^{M_{i}}_{W_{2}W_{3}})^{*}\right) \otimes M_{j}
\nno\\
&&
=\coprod_{j=1}^{k} \left(\coprod_{i=1}^{k} ({\cal
M}^{M_{j}}_{W_{1}M_{i}} \otimes {\cal
M}^{M_{i}}_{W_{2}W_{3}})^{*}\right) \otimes M_{j}
\end{eqnarray}
and
\begin{eqnarray}\label{(W1W2)W3}
\lefteqn{(W_{1}\boxtimes_{P(z)}W_{2})\boxtimes_{P(z)}W_{3}
=\left(\coprod_{i=1}^{k} M_{i} \otimes ({\cal
M}^{M_{i}}_{W_{1}W_{2}})^{*} \right) \boxtimes_{P(z)} W_{3}}
\nno\\
&&
=\coprod_{i=1}^{k} (M_{i} \boxtimes_{P(z)} W_{3}) \otimes ({\cal
M}^{M_{i}}_{W_{1}W_{2}})^{*}
\nno\\
&&
=\coprod_{i=1}^{k} \left(\coprod_{j=1}^{k} ({\cal
M}^{M_{j}}_{M_{i}W_{3}})^{*} \otimes M_{j}\right) \otimes ({\cal
M}^{M_{i}}_{W_{1}W_{2}})^{*}
\nno\\
&&
=\coprod_{j=1}^{k} \left(\coprod_{i=1}^{k} ({\cal
M}^{M_{j}}_{M_{i}W_{3}})^{*} \otimes ({\cal
M}^{M_{i}}_{W_{1}W_{2}})^{*}\right) \otimes M_{j}
\nno\\
&&
=\coprod_{j=1}^{k} \left(\coprod_{i=1}^{k} ({\cal
M}^{M_{j}}_{M_{i}W_{3}} \otimes {\cal
M}^{M_{i}}_{W_{1}W_{2}})^{*}\right) \otimes M_{j}.
\end{eqnarray}
These two $V$-modules will be equivalent if for each $j=1,\dots,k$,
their $M_j$-multiplicities are the same, that is, if
\begin{equation}\label{fusionrulerelation}
\sum_{i=1}^{k} N^{M_{j}}_{W_{1}M_{i}} N^{M_{i}}_{W_{2}W_{3}}=
\sum_{i=1}^{k} N^{M_{i}}_{W_{1}W_{2}} N^{M_{j}}_{M_{i}W_{3}}.
\end{equation}

{\it However,} knowing only that these two $V$-modules are equivalent
(knowing that $\boxtimes$ is ``associative'' in only a rough sense) is
far {}from enough.  What we need is a natural isomorphism between these
two modules analogous to the natural isomorphism
\begin{equation}\label{calWassociativity}
{\cal W}_1 \otimes ({\cal W}_2 \otimes {\cal W}_3)
\stackrel{\sim}{\longrightarrow} ({\cal W}_1 \otimes {\cal W}_2) \otimes
{\cal W}_3
\end{equation}
of vector spaces ${\cal W}_i$ determined by the natural condition
\begin{equation}\label{wassociativity}
w_{(1)} \otimes (w_{(2)} \otimes w_{(3)}) \mapsto (w_{(1)} \otimes w_{(2)})
\otimes w_{(3)}
\end{equation}
on elements (recall the Introduction).  Suppose that ${\cal W}_1,$
${\cal W}_2$ and ${\cal W}_3$ are finite-dimensional completely
reducible modules for some Lie algebra.  Then we of course have the
analogue of the relation (\ref{fusionrulerelation}).  But knowing the
equality of these multiplicities certainly does not give the natural
isomorphism (\ref{calWassociativity})--(\ref{wassociativity}).

Our intent to construct a natural isomorphism between the spaces
(\ref{W1(W2W3)}) and (\ref{(W1W2)W3}) (under suitable conditions) in
fact provides a guide to what we need to do.  In (\ref{W1(W2W3)}),
each space ${\cal M}^{M_{j}}_{W_{1}M_{i}} \otimes {\cal
M}^{M_{i}}_{W_{2}W_{3}}$ suggests combining an intertwining map ${\cal
Y}_1$ of type ${M_{j}}\choose {W_{1}M_{i}}$ with an intertwining map
${\cal Y}_2$ of type ${M_{i}}\choose {W_{2}W_{3}}$, presumably by
composition:
\begin{equation}\label{Y1zY2z}
{\cal Y}_1 (w_{(1)},z){\cal Y}_2 (w_{(2)},z).
\end{equation}
But this will not work, since this composition does not exist because
the relevant formal series in $z$ does not converge; we must instead
take
\begin{equation}\label{Y1z1Y2z2}
{\cal Y}_1 (w_{(1)},z_1){\cal Y}_2 (w_{(2)},z_2),
\end{equation}
where the complex numbers $z_1$ and $z_2$ are such that
\[
|z_1|>|z_2|>0,
\]
by analogy with, and generalizing, the situation in Corollary
\ref{dualitywithcovergence}.  The composition (\ref{Y1z1Y2z2}) must be
understood using convergence and ``matrix coefficients,'' again as in
Corollary \ref{dualitywithcovergence}.

Similarly, in (\ref{(W1W2)W3}), each space ${\cal
M}^{M_{j}}_{M_{i}W_{3}} \otimes {\cal M}^{M_{i}}_{W_{1}W_{2}}$
suggests combining an intertwining map ${\cal Y}^1$ of type
${M_{j}}\choose {M_{i}W_{3}}$ with an intertwining map of type ${\cal
Y}^2$ of type ${M_{i}}\choose {W_{1}W_{2}}$:
\[
{\cal Y}^1 ({\cal Y}^2 (w_{(1)},z_1-z_2)w_{(2)},z_2),
\]
a (convergent) iterate of intertwining maps as in 
(\ref{associativitywithz1,z2}), with
\[
|z_2|>|z_1-z_2|>0,
\]
{\it not}
\begin{equation}
{\cal Y}^1 ({\cal Y}^2 (w_{(1)},z)w_{(2)},z),
\end{equation}
which fails to converge.

The natural way to construct a natural associativity isomorphism
between (\ref{W1(W2W3)}) and (\ref{(W1W2)W3}) will in fact, then, be
to implement a correspondence of the type
\begin{equation}\label{YY=Y(Y)}
{\cal Y}_1 (w_{(1)},z_1){\cal Y}_2 (w_{(2)},z_2)={\cal Y}^1 ({\cal
Y}^2 (w_{(1)},z_1-z_2)w_{(2)},z_2),
\end{equation}
as we have previewed in the Introduction (formula (\ref{yyyy2})) and
also in (\ref{associativitywithz1,z2}).  Formula (\ref{YY=Y(Y)})
expresses the existence and associativity of the general
nonmeromorphic operator product expansion, as discussed in Remark
\ref{OPE}.  Note that this viewpoint shows that we should not try
directly to construct a natural isomorphism
\begin{equation}
W_1\boxtimes_{P(z)} (W_2\boxtimes_{P(z)} W_3)
\stackrel{\sim}{\longrightarrow} (W_1 \boxtimes_{P(z)}
W_2)\boxtimes_{P(z)} W_3,
\end{equation}
but rather a natural isomorphism
\begin{equation}\label{naturalassociso}
W_1\boxtimes_{P(z_1)} (W_2\boxtimes_{P(z_2)} W_3)
\stackrel{\sim}{\longrightarrow} (W_1 \boxtimes_{P(z_1-z_2)}
W_2)\boxtimes_{P(z_2)} W_3.
\end{equation}
This is what we will actually do in this work, in the general
logarithmic, not-necessarily-finitely-reductive case, under suitable
conditions.  The natural isomorphism (\ref{naturalassociso}) will act
as follows on elements of the completions of the relevant
(generalized) modules:
\begin{equation}
w_{(1)}\boxtimes_{P(z_1)} (w_{(2)}\boxtimes_{P(z_2)} w_{(3)})
\mapsto (w_{(1)} \boxtimes_{P(z_1-z_2)}
w_{(2)})\boxtimes_{P(z_2)} w_{(3)},
\end{equation}
implementing the strategy suggested by the classical natural
isomorphism (\ref{calWassociativity})--(\ref{wassociativity}).  Recall
that we previewed this strategy in the Introduction.

It turns out that in order to carry out this program, including the
construction of equalities of the type (\ref{YY=Y(Y)}) (the existence
and associativity of the nonmeromorphic operator product expansion) in
general, we cannot use the realization of the $P(z)$-tensor product
given in Proposition \ref{construcofPztensorprod-finredcase}, {\it
even when} $V$ {\it is a finitely reductive vertex operator algebra.}
As in \cite{tensor1}--\cite{tensor3} and \cite{tensor4}, what we do
instead is to construct $P(z)$-tensor products in a completely
different way (even in the finitely reductive case), a way that allows
us to also construct the natural associativity isomorphisms.  Section
5 is devoted to this construction of $P(z)$- (and $Q(z)$-)tensor
products.

\subsection{$Q(z)$-intertwining maps and the notion of $Q(z)$-tensor product}

We now generalize the notion of $Q(z)$-tensor product of modules {}from
\cite{tensor1} to the setting of the present work, parallel to what we
did for the $P(z)$-tensor product above.  Here we give only the
results that we will need later.  Other results similar to those for
$P(z)$-tensor products certainly also carry over to the case of
$Q(z)$, for example, the results above on the finitely reductive case,
as were presented in \cite{tensor1}.

\begin{defi}\label{im:qimdef}{\rm
Let $(W_1, Y_{1})$, $(W_2, Y_{2})$ and $(W_3, Y_{3})$ 
be generalized $V$-modules.  A {\it $Q(z)$-intertwining map of type
${W_3\choose W_1\,W_2}$} is a linear map 
\[
I: W_1\otimes W_2 \to
\overline{W}_3
\]
such that the following conditions are satisfied: 
the {\it grading
compatibility condition}: for $\beta, \gamma\in \tilde{A}$ and
$w_{(1)}\in W_{1}^{(\beta)}$, $w_{(2)}\in W_{2}^{(\gamma)}$,
\begin{equation}\label{grad-comp-qz}
I(w_{(1)}\otimes w_{(2)})\in \overline{W_{3}^{(\beta+\gamma)}};
\end{equation}
the
{\em lower truncation condition:} for any elements
$w_{(1)}\in W_1$, $w_{(2)}\in W_2$, and any $n\in {\mathbb C}$,
\begin{equation}\label{imq:ltc}
\pi_{n-m}I(w_{(1)}\otimes w_{(2)})=0\;\;\mbox{ for }\;m\in {\mathbb
N}\;\mbox{ sufficiently large}
\end{equation}
(which follows {}from (\ref{grad-comp-qz}), in view of the
grading restriction condition (\ref{set:dmltc}); cf. (\ref{im:ltc}));
the {\em Jacobi identity}:
\begin{eqnarray}\label{imq:def}
\lefteqn{z^{-1}\delta\left(\frac{x_1-x_0}{z}\right)
Y^o_3(v, x_0)I(w_{(1)}\otimes w_{(2)})}\nonumber\\
&&=x_0^{-1}\delta\left(\frac{x_1-z}{x_0}\right)
I(Y_1^{o}(v, x_1)w_{(1)}\otimes w_{(2)})\nonumber\\
&&\hspace{2em}-x_0^{-1}\delta\left(\frac{z-x_1}{-x_0}\right)
I(w_{(1)}\otimes Y_2(v, x_1)w_{(2)})
\end{eqnarray}
for $v\in V$, $w_{(1)}\in W_1$ and $w_{(2)}\in W_2$ (recall
(\ref{yo}) for the notation $Y^{o}$, and note that the 
left-hand side of (\ref{imq:def}) is
meaningful because any infinite linear combination of $v_n$ of the
form $\sum_{n<N}a_nv_n$ ($a_n\in {\mathbb C}$) acts on any
$I(w_{(1)}\otimes w_{(2)})$, in view of (\ref{imq:ltc})); and the {\em
${\mathfrak s}{\mathfrak l}(2)$-bracket relations}: for any
$w_{(1)}\in W_1$ and $w_{(2)}\in W_2$,
\begin{eqnarray}\label{imq:Lj}
L(-j)I(w_{(1)}\otimes w_{(2)})&=&\sum_{i=0}^{j+1}{j+1\choose
i}(-z)^iI((L(-j+i)w_{(1)})\otimes w_{(2)})\nno\\
&&-\sum_{i=0}^{j+1}{j+1\choose i}(-z)^iI(w_{(1)}\otimes
L(j-i)w_{(2)})
\end{eqnarray}
for $j=-1, 0$ and $1$ (note that if $V$ is in fact a conformal vertex
algebra, this follows automatically {}from (\ref{imq:def}) by setting
$v=\omega$ and taking $\res_{x_1}\res_{x_0}x_0^{j+1}$). The vector
space of $Q(z)$-intertwining maps of type ${W_3\choose W_1\,W_2}$ is
denoted by
\[
{\cal M}[Q(z)]^{W_3}_{W_1W_2}.
\]
}
\end{defi}

\begin{rema}\label{Q(z)geometry}{\rm
As was explained in \cite{tensor1}, the symbol $Q(z)$ represents the
Riemann sphere ${\mathbb C}\cup \{ \infty \}$ with one negatively
oriented puncture at $z$ and two ordered positively oriented punctures
at $\infty$ and $0$, with local coordinates $w-z$, $1/w$ and $w$,
respectively, vanishing at these punctures.  In fact, this structure
is conformally equivalent to the Riemann sphere ${\mathbb C}\cup \{
\infty \}$ with one negatively oriented puncture at $\infty$ and two
ordered positively oriented punctures $1/z$ and $0$, with local
coordinates $z/(zw-1)$, $(zw-1)/z^2w$ and $z^2w/(zw-1)$ vanishing at
$\infty$, $1/z$ and $0$, respectively. }
\end{rema}

\begin{rema}{\rm
In the case of $\C$-graded ordinary modules for a vertex operator
algebra, where the grading restriction condition (\ref{Wn+k=0}) for a
module $W$ is replaced by the (more restrictive) condition
\begin{equation}
W_{(n)}=0 \;\; \mbox { for }\;n\in {\C}\;\mbox{ with sufficiently
negative real part}
\end{equation}
as in \cite{tensor1} (and where, in our context, the abelian groups
$A$ and $\tilde{A}$ are trivial), the notion of $Q(z)$-intertwining
map above agrees with the earlier one introduced in \cite{tensor1}; in
this case, the conditions (\ref{grad-comp-qz}) and (\ref{imq:ltc}) are
automatic.}
\end{rema}

\begin{rema}{\rm (cf. Remark \ref{Pintwmaplowerbdd})  
If $W_3$ in Definition \ref{im:qimdef} is lower bounded, then
(\ref{imq:ltc}) can be strengthened to:
\begin{equation}\label{QpinI=0}
\pi_n I(w_{(1)}\otimes w_{(2)})=0\;\;\mbox{ for }
\;\Re{(n)}\;\mbox{ sufficiently negative.}
\end{equation}}
\end{rema}

In view of Remarks \ref{P(z)geometry} and \ref{Q(z)geometry}, we can
now give a natural correspondence between $P(z)$- and
$Q(z)$-intertwining maps.  (See the next three results.)  Recall that
since our generalized $V$-modules are strongly graded, we have
contragredient generalized modules of generalized modules.

\begin{propo}\label{qp:qp}
Let $I: W_1\otimes W_2\to \overline{W}_3$ and $J: W'_3\otimes W_2\to
\overline{W'_1}$ be linear maps related to each other by:
\begin{equation}\label{qz:qtop}
\langle w_{(1)},J(w'_{(3)}\otimes w_{(2)})\rangle=
\langle w'_{(3)},I(w_{(1)}\otimes w_{(2)})\rangle
\end{equation}
for any $w_{(1)}\in W_1$, $w_{(2)}\in W_2$ and $w'_{(3)}\in
W'_3$. Then $I$ is a $Q(z)$-intertwining map of type ${W_3\choose
W_1\,W_2}$ if and only if $J$ is a $P(z)$-intertwining map of type
${W'_1\choose W'_3\,W_2}$.
\end{propo}
\pf Suppose that $I$ is a $Q(z)$-intertwining map of type ${W_3\choose
W_1\,W_2}$.  We shall show that $J$ is a $P(z)$-intertwining map of
type ${W'_1\choose W'_3\,W_2}$.

Since $I$ satisfies the grading compatibility condition, it is clear
that $J$ also satisfies this condition.  For the lower truncation
condition for $J$, it suffices to show that for any $w_{(2)}\in
W_2^{(\beta)}$ and $w'_{(3)}\in (W'_3)^{(\gamma)}$, where $\beta, \gamma
\in \tilde{A}$, and any $n\in{\mathbb C}$,
\[
\langle \pi_{[n-m]}W_1^{(-\beta -\gamma)}, J(w'_{(3)}\otimes
w_{(2)})\rangle=0
\]
for $m\in{\mathbb N}$ sufficiently large, or that
\begin{equation}\label{qz:Jltrp}
\langle w'_{(3)},I(\pi_{[n-m]}W_1^{(-\beta -\gamma)} \otimes
w_{(2)})\rangle=0\;\;\mbox{ for }\;m\in{\mathbb
N}\;\;\mbox{sufficiently large.}
\end{equation}
But (\ref{qz:Jltrp}) follows immediately {}from (\ref{set:dmltc}).

Now we prove the Jacobi identity for $J$. The Jacobi identity for $I$
gives
\begin{eqnarray}\label{qz:jcba}
\lefteqn{z^{-1}\delta\left(\frac{x_1-x_0}{z}\right)\langle w'_{(3)},
Y^o_3(v, x_0)I(w_{(1)}\otimes w_{(2)})\rangle}\nno\\
&&=x_0^{-1}\delta\left(\frac{x_1-z}{x_0}\right)\langle w'_{(3)},
I(Y_1^{o}(v, x_1)w_{(1)}\otimes w_{(2)})\rangle\nno\\
&&\hspace{2em}-x_0^{-1}\delta\left(\frac{z-x_1}{-x_0} \right)\langle
w'_{(3)},I(w_{(1)}\otimes Y_2(v,x_1)w_{(2)})\rangle
\end{eqnarray}
for any $v\in V$, $w_{(1)}\in W_1$, $w_{(2)}\in W_2$ and $w'_{(3)}\in
W'_3$. By (\ref{y'}) the left-hand side is equal to
\[
z^{-1}\delta\left(\frac{x_1-x_0}{z}\right)\langle Y'_3(v,x_0)w'_{(3)},
I(w_{(1)}\otimes w_{(2)})\rangle
\]
So by (\ref{qz:qtop}), the identity (\ref{qz:jcba}) can be written as
\begin{eqnarray*}
\lefteqn{z^{-1}\delta\left(\frac{x_1-x_0}{z}\right)\langle w_{(1)},
J(Y'_3(v,x_0)w'_{(3)}\otimes w_{(2)})\rangle}\nno\\
&&=x_0^{-1}\delta\left(\frac{x_1-z}{x_0}\right)\langle Y_1^{o}(v,
x_1)w_{(1)}, J(w'_{(3)}\otimes w_{(2)})\rangle\nno\\
&&\hspace{2em}-x_0^{-1}\delta\left(\frac{z-x_1}{-x_0} \right)\langle
w_{(1)},J(w'_{(3)}\otimes Y_2(v,x_1)w_{(2)})\rangle.
\end{eqnarray*}
Applying (\ref{y'}) to the first term of the right-hand side
we see that this can be written as
\begin{eqnarray*}
\lefteqn{z^{-1}\delta\left(\frac{x_1-x_0}{z}\right)\langle w_{(1)},
J(Y'_3(v,x_0)w'_{(3)}\otimes w_{(2)})\rangle}\nno\\
&&=x_0^{-1}\delta\left(\frac{x_1-z}{x_0}\right)\langle w_{(1)},
Y'_1(v, x_1)J(w'_{(3)}\otimes w_{(2)})\rangle\nno\\
&&\hspace{2em}-x_0^{-1}\delta\left(\frac{z-x_1}{-x_0} \right)\langle
w_{(1)},J(w'_{(3)}\otimes Y_2(v,x_1)w_{(2)})\rangle
\end{eqnarray*}
for any $v\in V$, $w_{(1)}\in W_1$, $w_{(2)}\in W_2$ and $w'_{(3)}\in
W'_3$. This is exactly the Jacobi identity for $J$.

The ${\mathfrak s}{\mathfrak l}(2)$-bracket relations can be proved similarly,
as follows: The ${\mathfrak s}{\mathfrak l}(2)$-bracket relations for $I$ give
\begin{eqnarray*}
\langle w'_{(3)},L(-j)I(w_{(1)}\otimes w_{(2)})\rangle&=&
\sum_{i=0}^{j+1}{j+1\choose i}(-z)^i\langle
w'_{(3)},I((L(-j+i)w_{(1)})\otimes w_{(2)})\rangle\nno\\
&&-\sum_{i=0}^{j+1}{j+1\choose i}(-z)^i\langle
w'_{(3)},I(w_{(1)}\otimes L(j-i)w_{(2)})\rangle
\end{eqnarray*}
for any $w_{(1)}\in W_1$, $w_{(2)}\in W_2$, $w'_{(3)}\in W'_3$ and
$j=-1, 0, 1$. Using (\ref{L'(n)}) and then applying (\ref{qz:qtop}) we
get
\begin{eqnarray*}
\langle w_{(1)},J(L'(j)w'_{(3)}\otimes w_{(2)})\rangle&=&
\sum_{i=0}^{j+1}{j+1\choose i}(-z)^i\langle L(-j+i)w_{(1)},
J(w'_{(3)}\otimes w_{(2)})\rangle\nno\\
&&-\sum_{i=0}^{j+1}{j+1\choose i}(-z)^i\langle
w_{(1)},J(w'_{(3)}\otimes L(j-i)w_{(2)})\rangle,
\end{eqnarray*}
or
\begin{eqnarray*}
J(L'(j)w'_{(3)}\otimes w_{(2)})&=& \sum_{i=0}^{j+1}{j+1\choose
i}(-z)^i L(j-i)J(w'_{(3)}\otimes w_{(2)})\nno\\
&&-\sum_{i=0}^{j+1}{j+1\choose i}(-z)^iJ(w'_{(3)}\otimes
L(j-i)w_{(2)}),
\end{eqnarray*}
for $j=-1,0,1$. This is the alternative form (\ref{im:Lj2}) of the
${\mathfrak s} {\mathfrak l}(2)$-bracket relations for $J$. Hence $J$ is a
$P(z)$-intertwining map.

The other direction of the proposition is proved by simply reversing
the order of the arguments. \epfv

Let $W_{1}$, $W_{2}$ and $W_{3}$ be generalized $V$-modules, as above.
We shall call an element $\lambda$ of $(W_{1}\otimes W_{2}\otimes
W_{3})^{*}$ {\it $\tilde{A}$-compatible} if
\[
\lambda((W_{1})^{(\beta)}\otimes (W_{2})^{(\gamma)}\otimes
(W_{3})^{(\delta)})=0
\]
for $\beta, \gamma, \delta\in \tilde{A}$ satisfying
\[
\beta+\gamma+\delta\ne 0.
\]
Recall {}from Definitions \ref{Wbardef} and
\ref{defofWprime} that for a generalized $V$-module $W$,
$\overline{W'}$ can be viewed as a (usually proper) subspace of
$W^{*}$.  We shall call a linear map
\[
I: W_{1}\otimes W_{2} \rightarrow W_{3}^{*}
\]
{\it $\tilde{A}$-compatible} if its image lies in $\overline{W_{3}'}$,
that is,
\begin{equation}\label{IAtildecompat}
I:W_{1}\otimes W_{2} \rightarrow \overline{W_{3}'},
\end{equation}
and if $I$ satisfies the usual grading compatibility condition
(\ref{grad-comp}) or (\ref{grad-comp-qz}) for $P(z)$- or
$Q(z)$-intertwining maps.  
Now an element $\lambda$ of $(W_{1}\otimes
W_{2}\otimes W_{3})^{*}$ amounts exactly to a linear map
\[
I_{\lambda}:W_{1}\otimes W_{2} \rightarrow W_{3}^{*}.
\]
If $\lambda$ is $\tilde{A}$-compatible, then for $w_{(1)}\in
W_{1}^{(\beta)}$, $w_{(2)}\in W_{2}^{(\gamma)}$ and $w_{(3)}\in
W_{3}^{(\delta)}$ such that
\[
\delta\ne -(\beta +\gamma),
\]
we have
\[
\langle w_{(3)}, I_{\lambda}(w_{(1)}\otimes w_{(2)})\rangle=
\lambda(w_{(1)}\otimes w_{(2)}\otimes w_{(3)})=0,
\]
so that
\[
I_{\lambda}(w_{(1)}\otimes w_{(2)})\in \overline{(W_{3}')^{(\beta+\gamma)}}
\]
and $I_{\lambda}$ is
$\tilde{A}$-compatible.  Similarly, if $I_{\lambda}$ is
$\tilde{A}$-compatible, then so is $\lambda$.  Thus we have the
following straightforward result relating $\tilde{A}$-compatibility of
$\lambda$ with that of $I_{\lambda}$:

\begin{lemma}\label{4.36}
The linear functional $\lambda\in (W_{1}\otimes W_{2}\otimes
W_{3})^{*}$ is $\tilde{A}$-compatible if and only if $I_{\lambda}$ is
$\tilde{A}$-compatible.  The map given by $\lambda\mapsto I_{\lambda}$
is the unique linear isomorphism {}from the space of
$\tilde{A}$-compatible elements of $(W_{1}\otimes W_{2}\otimes
W_{3})^{*}$ to the space of $\tilde{A}$-compatible linear maps {}from
$W_{1}\otimes W_{2}$ to $\overline{W_{3}'}$ such that
\[
\langle w_{(3)}, I_{\lambda}(w_{(1)}\otimes w_{(2)})\rangle
=\lambda(w_{(1)}\otimes w_{(2)}\otimes w_{(3)})
\]
for $w_{(1)}\in W_{1}$, $w_{(2)}\in W_{2}$ and $w_{(3)}\in W_{3}$.
Similarly, there are canonical linear isomorphisms {}from the space of
$\tilde{A}$-compatible elements of $(W_{1}\otimes W_{2}\otimes
W_{3})^{*}$ to the space of $\tilde{A}$-compatible linear maps {}from
$W_{1}\otimes W_{3}$ to $\overline{W_{2}'}$ and to the space of
$\tilde{A}$-compatible linear maps {}from $W_{2}\otimes W_{3}$ to
$\overline{W_{1}'}$ satisfying the corresponding conditions.  In
particular, there is a canonical linear isomorphism {}from the space of
$\tilde{A}$-compatible linear maps {}from $W_{1}\otimes W_{2}$ to
$\overline{W_{3}}$ to the space of $\tilde{A}$-compatible linear maps
{}from $W_{3}'\otimes W_{2}$ to $\overline{W_{1}'}$ given by
(\ref{qz:qtop}).  \epf
\end{lemma}

Using this lemma and Proposition \ref{qp:qp}, we have:

\begin{corol}\label{Q(z)P(z)iso}
The formula (\ref{qz:qtop}) gives a canonical linear isomorphism
between the space of $Q(z)$-intertwining maps of type ${W_3\choose
W_1\,W_2}$ and the space of $P(z)$-intertwining maps of type
${W'_1\choose W'_3\,W_2}$.  \epf 
\end{corol}

\begin{rema}{\rm If the generalized modules under consideration are
lower bounded, then the spaces of intertwining maps satisfy the
stronger conditions (\ref{PpinI=0}) and (\ref{QpinI=0}).}
\end{rema}

We can now use Proposition \ref{im:correspond} together with
Proposition \ref{qp:qp} and Corollary \ref{Q(z)P(z)iso} to construct a
correspondence between the logarithmic intertwining operators of type
${W'_1}\choose {W'_3W_2}$ and the $Q(z)$-intertwining maps of type
${W_3}\choose {W_1W_2}$; this generalizes the corresponding result in
the finitely reductive case, with ordinary modules, in \cite{tensor1}.
Fix an integer $p$.  Let ${\cal Y}$ be a logarithmic intertwining
operator of type ${W'_1}\choose {W'_3W_2}$, and use (\ref{log:IYp}) to
define a linear map
\[
I_{{\cal Y}, p}: W'_3\otimes W_2\to \overline{W_{1}'};
\]
by Proposition \ref{im:correspond}, this is a
$P(z)$-intertwining map of the same type.  Then use Proposition
\ref{qp:qp} and Corollary \ref{Q(z)P(z)iso} to define a
$Q(z)$-intertwining map
\[
I^{Q(z)}_{{\cal Y}, p}: W_1\otimes W_2\to
\overline{W}_3
\]
of type ${W_3}\choose {W_1W_2}$ (uniquely) by
\begin{eqnarray}\label{imq:IYp}
\lefteqn{\langle w'_{(3)}, I^{Q(z)}_{{\cal Y}, p}(w_{(1)}\otimes
w_{(2)})\rangle_{W_3} = \langle w_{(1)}, I_{{\cal Y},
p}(w'_{(3)}\otimes w_{(2)})\rangle_{W'_1}}\nno\\
&&=\langle w_{(1)}, {\cal Y}(w'_{(3)},
e^{l_{p}(z)})w_{(2)}\rangle_{W'_1}
\;\;\;\;\;\;\;\;\;\;\;\;\;\;\;\;\;\;\;\;\;\;\;\;\;\;\;\;\;\;\;\;\;\;\;\;\;\;\;
\end{eqnarray}
for all $w_{(1)}\in W_1$, $w_{(2)}\in W_2$, $w'_{(3)}\in W'_3$. (We
are using the symbol $Q(z)$ to distinguish this {}from the $P(z)$ case
above.)  Then the correspondence
\[
{\cal Y} \mapsto I^{Q(z)}_{{\cal Y},p}
\]
is an isomorphism {}from ${\cal V}^{W_{1}'}_{W_{3}'W_{2}}$ to ${\cal
M}[Q(z)]^{W_3}_{W_1W_2}$.  {}From Proposition \ref{im:correspond} and
(\ref{recover}), its inverse is given by sending a $Q(z)$-intertwining
map $I$ of type ${W_3}\choose {W_1W_2}$ to the logarithmic
intertwining operator
\[
{\cal Y}^{Q(z)}_{I,p}:W'_3\otimes W_2\to W'_1[\log x]\{x\}
\]
defined by
\begin{eqnarray*}
\lefteqn{\langle w_{(1)}, {\cal Y}^{Q(z)}_{I,p}(w'_{(3)},
x)w_{(2)}\rangle_{W'_1}}\\ 
&&=\langle y^{-L'(0)}x^{-L'(0)}w'_{(3)},
I(y^{L(0)}x^{L(0)}w_{(1)}\otimes
y^{-L(0)}x^{-L(0)}w_{(2)})\rangle_{W_3}\lbar_{y=e^{-l_{p}(z)}}
\end{eqnarray*}
for any $w_{(1)}\in W_1$, $w_{(2)}\in W_2$ and $w'_{(3)}\in
W'_3$.  Thus we have:

\begin{propo}\label{Q-cor}
For $p\in {\mathbb Z}$, the correspondence
\[
{\cal Y}\mapsto I^{Q(z)}_{{\cal Y}, p}
\]
is a linear isomorphism {}from the space ${\cal
V}^{W'_1}_{W'_3W_2}$ of logarithmic intertwining operators of type
${W'_1}\choose {W'_3\; W_2}$ to the space ${\cal
M}[Q(z)]^{W_3}_{W_1W_2}$ of $Q(z)$-intertwining maps of type
${W_3}\choose {W_1W_2}$. Its inverse is given by
\[
I\mapsto {\cal Y}^{Q(z)}_{I, p}.
\]
\epf
\end{propo}

\begin{rema}{\rm If the generalized modules under consideration are
lower bounded, then the stronger conditions (\ref{repartbounded}) and
(\ref{QpinI=0}) hold.}
\end{rema}

We now give the definition of $Q(z)$-tensor product.

\begin{defi}\label{qz-product}{\rm
Let ${\cal C}_{1}$ be either ${\cal M}_{sg}$ or ${\cal GM}_{sg}$.
For $W_1, W_2\in \ob{\cal C}_{1}$, a {\it $Q(z)$-product of $W_1$ and
$W_2$} is an object $(W_3, Y_3)$ of ${\cal C}_{1}$ together with a
$Q(z)$-intertwining map $I_3$ of type ${W_3}\choose {W_1W_2}$. We
denote it by $(W_3, Y_3; I_3)$ or simply by $(W_3, I_3)$.  Let
$(W_4,Y_4; I_4)$ be another $Q(z)$-product of $W_1$ and $W_2$. A {\em
morphism} {}from $(W_3, Y_3; I_3)$ to $(W_4, Y_4; I_4)$ is a module map
$\eta$ {}from $W_3$ to $W_4$ such that the
diagram
\begin{center}
\begin{picture}(100,60)
\put(-5,0){$\overline W_3$}
\put(13,4){\vector(1,0){104}}
\put(119,0){$\overline W_4$}
\put(41,50){$W_1\otimes W_2$}
\put(61,45){\vector(-3,-2){50}}
\put(68,45){\vector(3,-2){50}}
\put(65,8){$\bar\eta$}
\put(20,27){$I_3$}
\put(98,27){$I_4$}
\end{picture}
\end{center}
commutes, that is,
\[
I_4=\overline{\eta}\circ I_3.
\]
where, as before, $\overline{\eta}$ is the natural map {from}
$\overline{W}_3$ to $\overline{W}_4$ extending $\eta$. }
\end{defi}

\begin{defi}\label{qz-tp}
{\rm Let ${\cal C}$ be a full subcategory of either ${\cal M}_{sg}$ or
${\cal GM}_{sg}$. For $W_1, W_2\in \ob{\cal C}$, a {\em $Q(z)$-tensor
product of $W_1$ and $W_2$ in ${\cal C}$} is a $Q(z)$-product $(W_0,
Y_0; I_0)$ with $W_0\in{\rm ob\,}{\cal C}$ such that for any
$Q(z)$-product $(W,Y;I)$ with $W\in{\rm ob\,}{\cal C}$, there is a
unique morphism {}from $(W_0, Y_0; I_0)$ to $(W,Y;I)$. Clearly, a
$Q(z)$-tensor product of $W_1$ and $W_2$ in ${\cal C}$, if it exists,
is unique up to unique isomorphism.  In this case we will denote it by
\[
(W_1\boxtimes_{Q(z)} W_2, Y_{Q(z)}; \boxtimes_{Q(z)})
\]
and call the object
\[
(W_1\boxtimes_{Q(z)} W_2, Y_{Q(z)})
\]
the {\em $Q(z)$-tensor
product (generalized) module of $W_1$ and $W_2$ in ${\cal C}$}. Again
we will skip the phrase ``in ${\cal C}$'' if the category ${\cal C}$
under consideration is clear in context. }
\end{defi}

The following immediate consequence of Definition \ref{qz-tp}
and Proposition \ref{Q-cor}
relates module maps {}from a
$Q(z)$-tensor product module with $Q(z)$-intertwining maps and
logarithmic intertwining operators:

\begin{propo}
Suppose that $W_1\boxtimes_{Q(z)}W_2$ exists. We have a natural
isomorphism
\begin{eqnarray*}
\hom_{V}(W_1\boxtimes_{Q(z)}W_2, W_3)&\stackrel{\sim}{\to}&
{\cal M}[Q(z)]^{W_3}_{W_1W_2}\\ \eta&\mapsto& \overline{\eta}\circ
\boxtimes_{Q(z)}
\end{eqnarray*}
and for $p\in {\mathbb Z}$, a natural isomorphism
\begin{eqnarray*}
\hom_{V}(W_1\boxtimes_{Q(z)} W_2, W_3)&
\stackrel{\sim}{\rightarrow}& {\cal V}^{W'_1}_{W'_3W_2}\\ \eta&\mapsto
& {\cal Y}^{Q(z)}_{\eta, p}
\end{eqnarray*}
where ${\cal Y}^{Q(z)}_{\eta, p}={\cal Y}^{Q(z)}_{I, p}$ with
$I=\overline{\eta}\circ \boxtimes_{Q(z)}$.\epf
\end{propo}

Suppose that the $Q(z)$-tensor product $(W_1\boxtimes_{Q(z)} W_2,
Y_{Q(z)}; \boxtimes_{Q(z)})$ of $W_1$ and $W_2$ exists.  We will
sometimes denote the action of the canonical $Q(z)$-intertwining map
\begin{equation}\label{q-actionofboxtensormap}
w_{(1)} \otimes w_{(2)}\mapsto\boxtimes_{Q(z)}(w_{(1)} \otimes
w_{(2)})=\boxtimes_{Q(z)}(w_{(1)}, z)w_{(2)} \in 
\overline{W_1\boxtimes_{Q(z)}W_2}
\end{equation}
on elements simply by $w_{(1)}\boxtimes_{Q(z)}
w_{(2)}$:
\begin{equation}\label{q-boxtensorofelements}
w_{(1)}\boxtimes_{Q(z)} w_{(2)}=\boxtimes_{Q(z)}(w_{(1)} \otimes
w_{(2)})=\boxtimes_{Q(z)}(w_{(1)}, z)w_{(2)}.
\end{equation}

Using Propositions \ref{log:omega} and \ref{log:A}, we have the
following result, generalizing Proposition 4.9 and Corollary 4.10 in
\cite{tensor1}:

\begin{propo}\label{b-r}
For any integer $r$, there is a natural isomorphism
\[
B_{r}: {\cal V}^{W_3}_{W_1W_2}\to {\cal V}^{W'_1}_{W'_3W_2}
\]
defined by the condition that for any logarithmic intertwining
operator ${\cal Y}$ in ${\cal V}^{W_3}_{W_1W_2}$ and $w_{(1)}\in W_1$,
$w_{(2)}\in W_2$, $w'_{(3)}\in W'_3$,
\begin{eqnarray}\label{4.31}
\lefteqn{\langle w_{(1)}, B_{r}({\cal Y})(w'_{(3)}, x)
w_{(2)}\rangle_{W'_1}}\nno\\
&&=\langle e^{-x^{-1}L(1)}w'_{(3)}, {\cal Y}(e^{xL(1)}w_{(1)},
x^{-1})e^{-xL(1)}e^{(2r+1)\pi iL(0)}
(x^{-L(0)})^{2}w_{(2)}\rangle_{W_3}.
\end{eqnarray}
\end{propo}
\pf
{From} Proposition \ref{log:omega},  for any integer $r_{1}$ we have an
isomorphism $\Omega_{r_{1}}$ {from} ${\cal
V}^{W_{3}}_{W_{1}W_{2}}$ to ${\cal V}^{W_{3}}_{W_{2}W_{1}}$, and {from}
Proposition \ref{log:A},  for any integer $r_{2}$ we have an
isomorphism $A_{r_{2}}$ {from} ${\cal V}^{W_{3}}_{W_{2}W_{1}}$
to ${\cal V}^{W'_{1}}_{W_{2}W'_{3}}$. By Proposition \ref{log:omega} again, 
for any integer $r_{3}$ there is an  isomorphism, which we again denote 
$\Omega_{r_{3}}$, {from} ${\cal V}^{W'_{1}}_{W_{2}W'_{3}}$ to ${\cal
V}^{W'_{1}}_{W'_{3}W_{2}}$. Thus for any triple $(r_{1}, r_{2},
r_{3})$ of integers, we have an isomorphism $\Omega_{r_{3}}\circ A_{r_{2}}\circ
\Omega_{r_{1}}$ {from} ${\cal V}^{W_{3}}_{W_{1}W_{2}}$ to ${\cal
V}^{W'_{1}}_{W'_{3}W_{2}}$.  Let ${\cal Y}$ be a logarithmic intertwining
operator in ${\cal V}^{W_{3}}_{W_{1}W_{2}}$ and $w_{(1)}$, $w_{(2)}$,
$w'_{(3)}$ elements of $W_{1}$, $W_{2}$, $W'_{3}$, respectively.
{From} the definitions of $\Omega_{r_{1}}$, $A_{r_{2}}$ and
$\Omega_{r_{3}}$, we have
\begin{eqnarray}\label{7.29}
\lefteqn{\langle (\Omega_{r_{3}}\circ A_{r_{2}}\circ 
\Omega_{r_{1}})({\cal Y})(w'_{(3)}, x)w_{(2)}, w_{(1)}\rangle_{W_1}=}\nno\\
&&=\langle e^{xL(-1)}A_{r_{2}}(\Omega_{r_{1}}
({\cal Y}))(w_{(2)}, e^{(2r_{3}+1)\pi i}x)w'_{(3)}, 
w_{(1)}\rangle_{W_1}\nno\\
&&=\langle A_{r_{2}}(\Omega_{r_{1}}
({\cal Y}))(w_{(2)}, e^{(2r_{3}+1)\pi i}x)w'_{(3)}, 
e^{xL(1)}w_{(1)}\rangle_{W_1}\nno\\
&&=\langle w'_{(3)}, \Omega_{r_{1}}({\cal Y})(e^{-xL(1)}e^{(2r_{2}+1)\pi iL(0)}
e^{-2(2r_{3}+1)\pi iL(0)}(x^{-L(0)})^{2}w_{(2)}, \nno\\
&&\hspace{10em}e^{-(2r_{3}+1)\pi i}x^{-1})
e^{xL(1)}w_{(1)}\rangle_{W_3}
\nno\\
&&=\langle w'_{(3)}, e^{-x^{-1}L(-1)}{\cal Y}(e^{xL(1)}w_{(1)}, 
e^{(2r_{1}+1)\pi i}e^{-(2r_{3}+1)\pi i}x^{-1})\cdot \nno\\
&&\hspace{4em}\cdot e^{-xL(1)}e^{(2r_{2}+1)\pi iL(0)}
e^{-2(2r_{3}+1)\pi iL(0)}(x^{-L(0)})^{2}w_{(2)}\rangle_{W_3}
\nno\\
&&=\langle e^{-x^{-1}L(1)}w'_{(3)}, {\cal Y}(e^{xL(1)}w_{(1)}, 
e^{2(r_{1}-r_{3})\pi i}x^{-1})\cdot \nno\\
&&\hspace{6em}\cdot e^{-xL(1)}e^{(2(r_{2}-2r_{3}-1)+1)\pi iL(0)}
(x^{-L(0)})^{2}w_{(2)}\rangle_{W_3}.
\end{eqnarray}
{From} (\ref{7.29}) we see that $\Omega_{r_{3}}\circ A_{r_{2}}\circ
\Omega_{r_{1}}$ depends only on $r_{2}-2r_{3}-1$ and $r_{1}-r_{3}$,
and the operators $\Omega_{r_{3}}\circ A_{r_{2}}\circ \Omega_{r_{1}}$
with different $r_{1}-r_{3}$ but the same $r_{2}-2r_{3}-1$ differ
{from} each other only by automorphisms of ${\cal
V}^{W_{3}}_{W_{1}W_{2}}$ (recall Remarks \ref{log:fcf},
\ref{exponentialaVhom} and \ref{Ys1s2s3}).  Thus for our purpose, we
need only consider those isomorphisms such that $r_{1}-r_{3}=0$.
Given any integer $r$, we choose two integers $r_{2}$ and $r_{3}$ such
that $r=r_{2}-2r_{3}-1$ and we define
\begin{equation}
B_{r}=\Omega_{r_{3}}\circ A_{r_{2}}\circ \Omega_{r_{3}}.
\end{equation}
{From} (\ref{7.29}) we see that $B_{r}$ is independent of the choices of $r_{2}$ and
$r_{3}$ and that 
(\ref{4.31}) holds.
\epfv

Combining the last two results, we obtain:

\begin{corol}
 For any $W_1, W_2, W_3\in \ob{\cal C}$ such that
$W_1\boxtimes_{Q(z)}W_2$ exists and any integers $p$ and $r$, we have
a natural isomorphism
\begin{eqnarray}
\hom_{V}(W_1\boxtimes_{Q(z)} W_2, W_3)&
\stackrel{\sim}{\rightarrow}&
{\cal V}^{W_3}_{W_1W_2}\nno\\
\eta&\mapsto &B^{-1}_{r}({\cal Y}^{Q(z)}_{\eta, p}).\hspace{2em}\square
\end{eqnarray}
\end{corol}

\subsection{$P(z)$-tensor products and $Q(z^{-1})$-tensor products}

Here we prove the following result:

\begin{theo}\label{pz-qz-1}
Let $W_{1}$ and $W_{2}$ be objects of a full subcategory $\mathcal{C}$
of either $\mathcal{M}_{sg}$ or $\mathcal{GM}_{sg}$. Then the
$P(z)$-tensor product of $W_{1}$ and $W_{2}$ exists if and only if the
$Q(z^{-1})$-tensor product of $W_{1}$ and $W_{2}$ exists.
\end{theo}
\pf 
Recalling our choice of branch (\ref{branch1}), let
\[
p=-\frac{\log (z^{-1})+\log z}{2\pi i}.
\]
Then $p$ is an integer and we have
\[
-(\log (z^{-1})+2\pi pi)=\log z,
\]
and 
\[
e^{-n(\log (z^{-1})+2\pi pi)}=e^{n\log z}
\]
for $n\in \C$.

{}From Propositions \ref{im:correspond}, \ref{b-r} and
\ref{Q-cor}, we see that for $W_{1}, W_{2}, W_{3}\in \ob \mathcal{C}$,
there is a linear isomorphism $\mu_{W_{1}W_{2}}^{W_{3}}:
\mathcal{M}[P(z)]_{W_{1}W_{2}}^{W_{3}} \to
\mathcal{M}[Q(z^{-1})]_{W_{1}W_{2}}^{W_{3}}$ defined by
\[
\mu_{W_{1}W_{2}}^{W_{3}}(I)=I^{Q(z^{-1})}_{B_{2p}(\mathcal{Y}_{I,0}),p}
\]
for $I\in \mathcal{M}[P(z)]_{W_{1}W_{2}}^{W_{3}}$.  
By definition, $\mu_{W_{1}W_{2}}^{W_{3}}(I)$ is determined uniquely
by (recalling (\ref{branch1})--(\ref{branch2}))
\begin{eqnarray}\label{mu}
\lefteqn{\langle w_{(3)}', \mu_{W_{1}W_{2}}^{W_{3}}(I)(w_{(1)}
\otimes w_{(2)})\rangle}\nn
&&=\langle w_{(3)}', I^{Q(z^{-1})}_{B_{2p}(\mathcal{Y}_{I, 0}),p}
(w_{(1)}\otimes w_{(2)})\rangle\nn
&&=\langle w_{(1)}, B_{2p}(\mathcal{Y}_{I, 0})(w'_{(3)},
e^{l_p(z^{-1})})w_{(2)}\rangle\nn
&&=\langle w_{(1)}, B_{2p}(\mathcal{Y}_{I, 0})(w'_{(3)},
e^{\log (z^{-1})+2\pi pi})w_{(2)}\rangle\nn
&&=\langle e^{-zL(1)}w'_{(3)}, 
\mathcal{Y}_{I, 0}(e^{z^{-1}L(1)}w_{(1)}, e^{\log z})e^{-z^{-1}L(1)}
e^{(2(2p)+1)i\pi L(0)}e^{-2(\log z^{-1}+2\pi pi))L(0)}w_{(2)}\rangle\nn
&&=\langle e^{-zL(1)}w'_{(3)}, 
\mathcal{Y}_{I, 0}(e^{z^{-1}L(1)}w_{(1)}, z)e^{-z^{-1}L(1)}
e^{i\pi L(0)}e^{-2(\log z^{-1})L(0)}w_{(2)}\rangle\nn
&&=\langle e^{-zL(1)}w'_{(3)}, I((e^{z^{-1}L(1)}w_{(1)})
\otimes (e^{-z^{-1}L(1)}
e^{i\pi L(0)}e^{-2(\log z^{-1})L(0)}w_{(2)}))\rangle\nn
\end{eqnarray}
for $w_{(1)}\in W_{1}$, $w_{(2)}\in W_{2}$ and $w_{(3)}'\in W_{3}'$.
{}From (\ref{mu}), we also see that
for $J\in \mathcal{M}[Q(z^{-1})]_{W_{1}W_{2}}^{W_{3}}$,
$(\mu_{W_{1}W_{2}}^{W_{3}})^{-1}(J)$ is determined uniquely
by 
\begin{eqnarray}\label{mu-1}
\lefteqn{\langle w_{(3)}', (\mu_{W_{1}W_{2}}^{W_{3}})^{-1}(J)(w_{(1)}
\otimes w_{(2)})\rangle}\nn
&&=\langle e^{zL(1)}w'_{(3)}, 
J((e^{-z^{-1}L(1)}w_{(1)})
\otimes (e^{2(\log z^{-1})L(0)}
e^{-i\pi L(0)}e^{z^{-1}L(1)}w_{(2)}))\rangle\nn
\end{eqnarray}
for $w_{(1)}\in W_{1}$, $w_{(2)}\in W_{2}$ and $w_{(3)}'\in W_{3}'$.

Assume that the
$P(z)$-tensor product $(W_{1}\boxtimes_{P(z)} W_{2}, Y_{P(z)};
\boxtimes_{P(z)})$ exists.  Then
\[
\boxtimes_{Q(z^{-1})}=\mu_{W_{1}W_{2}}^{W_{1}\boxtimes_{P(z)}
W_{2}}(\boxtimes_{P(z)})=
I^{Q(z^{-1})}_{B_{2p}(\mathcal{Y}_{\boxtimes_{P(z)}, 0}),p}
\]
is a $Q(z^{-1})$-intertwining
map
of type ${W_{1}\boxtimes_{P(z)}
W_{2}\choose W_{1}W_{2}}$. 
We claim that $(W_{1}\boxtimes_{P(z)} W_{2}, Y_{P(z)};
\boxtimes_{Q(z^{-1})})$ is the $Q(z^{-1})$-tensor product of $W_{1}$ and
$W_{2}$.  

In fact, for any $Q(z^{-1})$-product $(W, Y; I)$ of $W_{1}$ and
$W_{2}$,
\[
(\mu_{W_{1}W_{2}}^{W})^{-1}(I)
=I_{B^{-1}_{2p}(\mathcal{Y}_{I,p}^{Q(z^{-1})}), 0}
\]
is a $P(z)$-intertwining map of type ${W\choose W_{1}W_{2}}$ and thus
$(W, Y; (\mu_{W_{1}W_{2}}^{W})^{-1}(I))$ is a $P(z)$-product of
$W_{1}$ and $W_{2}$. Since $(W_{1}\boxtimes_{P(z)} W_{2}, Y_{P(z)};
\boxtimes_{P(z)})$ is the $P(z)$-tensor product of $W_{1}$ and
$W_{2}$, there is a unique morphism of $P(z)$-products {}from
$(W_{1}\boxtimes_{P(z)} W_{2}, Y_{P(z)}; \boxtimes_{P(z)})$ to $(W, Y;
(\mu_{W_{1}W_{2}}^{W})^{-1}(I))$, that is, there exists a unique
module map 
\[
\eta^{P(z)}: W_{1}\boxtimes_{P(z)} W_{2} \to W
\]
such that
\[
(\mu_{W_{1}W_{2}}^{W})^{-1}(I)=\overline{\eta^{P(z)}}
\circ \boxtimes_{P(z)},
\]
or equivalently,
\begin{eqnarray}\label{pzt-qzt-equiv-1}
I&=&\mu_{W_{1}W_{2}}^{W}(\overline{\eta^{P(z)}}
\circ \boxtimes_{P(z)})\nn
&=&\mu_{W_{1}W_{2}}^{W}(\overline{\eta^{P(z)}}\circ 
(\mu_{W_{1}W_{2}}^{W_{1}\boxtimes_{P(z)}
W_{2}})^{-1}(\mu_{W_{1}W_{2}}^{W_{1}\boxtimes_{P(z)}
W_{2}}(\boxtimes_{P(z)})))\nn
&=&\mu_{W_{1}W_{2}}^{W}(\overline{\eta^{P(z)}}\circ 
(\mu_{W_{1}W_{2}}^{W_{1}\boxtimes_{P(z)}
W_{2}})^{-1}(\boxtimes_{Q(z^{-1})})).
\end{eqnarray}

{}From (\ref{mu}) and (\ref{mu-1}), we see that the right-hand side
of (\ref{pzt-qzt-equiv-1}) is determined uniquely by 
\begin{eqnarray}\label{pzt-qzt-equiv-2}
\lefteqn{\langle w', (\mu_{W_{1}W_{2}}^{W}(\overline{\eta^{P(z)}}\circ 
(\mu_{W_{1}W_{2}}^{W_{1}\boxtimes_{P(z)}
W_{2}})^{-1}(\boxtimes_{Q(z^{-1})})))(w_{(1)}\otimes w_{(2)})\rangle}\nn
&&=\langle e^{-zL(1)}w', (\overline{\eta^{P(z)}}\circ 
(\mu_{W_{1}W_{2}}^{W_{1}\boxtimes_{P(z)}
W_{2}})^{-1}(\boxtimes_{Q(z^{-1})}))\nn
&&\quad\quad\quad\quad\quad\quad\quad((e^{z^{-1}L(1)}w_{(1)})
\otimes (e^{-z^{-1}L(1)}
e^{i\pi L(0)}e^{-2(\log z^{-1})L(0)}w_{(2)}))\rangle\nn
&&=\langle e^{-zL(1)}w', \overline{\eta^{P(z)}}
((\mu_{W_{1}W_{2}}^{W_{1}\boxtimes_{P(z)}
W_{2}})^{-1}(\boxtimes_{Q(z^{-1})})\nn
&&\quad\quad\quad\quad\quad\quad\quad((e^{z^{-1}L(1)}w_{(1)})
\otimes (e^{-z^{-1}L(1)}
e^{i\pi L(0)}e^{-2(\log z^{-1})L(0)}w_{(2)})))\rangle\nn
&&=\langle (\eta^{P(z)})'(e^{-zL(1)}w'), 
(\mu_{W_{1}W_{2}}^{W_{1}\boxtimes_{P(z)}
W_{2}})^{-1}(\boxtimes_{Q(z^{-1})})\nn
&&\quad\quad\quad\quad\quad\quad\quad((e^{z^{-1}L(1)}w_{(1)})
\otimes (e^{-z^{-1}L(1)}
e^{i\pi L(0)}e^{-2(\log z^{-1})L(0)}w_{(2)}))\rangle\nn
&&=\langle (\eta^{P(z)})'(w'), 
\boxtimes_{Q(z^{-1})}(w_{(1)}
\otimes w_{(2)})\rangle\nn
&&=\langle w',  (\overline{\eta^{P(z)}}\circ \boxtimes_{Q(z^{-1})})
(w_{(1)}\otimes w_{(2)})\rangle
\end{eqnarray}
for $w_{(1)}\in W_{1}$, $w_{(2)}\in W_{2}$ and $w'\in W'$. {}From 
(\ref{pzt-qzt-equiv-1}) and (\ref{pzt-qzt-equiv-2}), we
see that 
\begin{equation}\label{pzt-qzt-equiv-3}
I=\overline{\eta^{P(z)}}\circ \boxtimes_{Q(z^{-1})}.
\end{equation}

We also need to show the uniqueness---that any module map $\eta:
W_{1}\boxtimes_{P(z)} W_{2}\to W$ such that $I=\overline{\eta}\circ
\boxtimes_{Q(z^{-1})}$ must be equal to $\eta^{P(z)}$.  For this,
it is sufficient to show that $\eta_1 = 0$, where
\[
\eta_1 = \eta^{P(z)} - \eta,
\]
given that 
\[
\overline{\eta_1}(w_{(1)}\boxtimes_{Q(z^{-1})}w_{(2)})=0
\]
for $w_{(1)}\in W_{1}$ and $w_{(2)}\in W_{2}$.  But for $w'\in
(W_{1}\boxtimes_{P(z)}W_{2})'$
\[
\langle e^{zL(1)} w', \overline{\eta_1}
(w_{(1)}\boxtimes_{Q(z^{-1})}w_{(2)})\rangle=0,
\]
so that
\[
\langle e^{zL(1)} {\eta_1}'(w'), 
w_{(1)}\boxtimes_{Q(z^{-1})}w_{(2)}\rangle
=\langle {\eta_1}' (e^{zL(1)} w'), 
w_{(1)}\boxtimes_{Q(z^{-1})}w_{(2)}\rangle = 0.
\]
{}From the definition of $\boxtimes_{Q(z)}$ and (\ref{mu}), we have
\begin{eqnarray*}
\lefteqn{\langle 
e^{zL(1)} {\eta_1}' (w'), w_{(1)}\boxtimes_{Q(z^{-1})}w_{(2)}\rangle}\nn
&&=\langle {\eta_1}' (w'), 
(e^{z^{-1}L(1)}w_{(1)})\boxtimes_{P(z)}(e^{-z^{-1}L(1)}
e^{i\pi L(0)}e^{-2(\log z^{-1})L(0)}w_{(2)})\rangle, \nn
\end{eqnarray*}
and thus
\begin{equation}\label{span2}
\langle {\eta_1}'(w'), 
(e^{z^{-1}L(1)}w_{(1)})\boxtimes_{P(z)}(e^{-z^{-1}L(1)}
e^{i\pi L(0)}e^{-2(\log z^{-1})L(0)}w_{(2)})\rangle=0.
\end{equation}
Since  $e^{z^{-1}L(1)}$ and $e^{-z^{-1}L(1)}
e^{i\pi L(0)}e^{-2(\log z^{-1})L(0)}$ are invertible operators
on $W_{1}$ and $W_{2}$, (\ref{span2}) for all $w_{(1)}\in W_{1}$, 
$w_{(2)}\in W_{2}$ is equivalent to 
\[
\langle {\eta_1}'(w'), w_{(1)}\boxtimes_{P(z)}w_{(2)}\rangle=0
\]
for all $w_{(1)}\in W_{1}$, $w_{(2)}\in W_{2}$.  Thus by
Proposition \ref{span},
\[
{\eta_1}'(w')=0
\]
for all homogeneous $w'$ and hence for all $w'$, showing that indeed
$\eta_1 = 0$ and proving the uniqueness of $\eta$.  Thus
$(W_{1}\boxtimes_{P(z)} W_{2}, Y_{P(z)}; \boxtimes_{Q(z^{-1})})$ is
the $Q(z^{-1})$-tensor product of $W_{1}$ and $W_{2}$.

Conversely, by essentially reversing these arguments we see that if
the $Q(z^{-1})$-tensor product of $W_{1}$ and $W_{2}$ exists, then so
does the $P(z)$-tensor product.  \epfv

{}From Theorem \ref{pz-qz-1} and Proposition \ref{4.19}, we
immediately obtain:

\begin{corol}\label{pz-qz}
Let $W_{1}$ and $W_{2}$ be objects of a full subcategory $\mathcal{C}$
of either $\mathcal{M}_{sg}$ or $\mathcal{GM}_{sg}$. Then the
$P(z)$-tensor product of $W_{1}$ and $W_{2}$ exists if and only if the
$Q(z)$-tensor product of $W_{1}$ and $W_{2}$ exists. \epf
\end{corol}

\begin{rema}{\rm
{}From the proof we see that as generalized $V$-modules,
$W_{1}\boxtimes_{P(z)} W_{2}$ and $W_{1}\boxtimes_{Q(z^{-1})} W_{2}$
are equivalent, but the main issue is that the intertwining maps
$\boxtimes_{P(z)}$ and $\boxtimes_{Q(z^{-1})}$, which encode the
geometric information, are very different; as generalized $V$-modules
{\it only}, $W_{1}\boxtimes_{P(z)} W_{2}$ and $W_{1}\boxtimes_{Q(z)}
W_{2}$ are equivalent.  Compare this with Remark
\ref{intwmapdependsongeomdata}.}
\end{rema}

%\section{References}

\bigskip

\noindent {\small \sc Department of Mathematics, Rutgers University,
Piscataway, NJ 08854 (permanent address)}

\noindent {\it and}

\noindent {\small \sc Beijing International Center for Mathematical Research,
Peking University, Beijing, China}

\noindent {\em E-mail address}: yzhuang@math.rutgers.edu

\vspace{1em}

\noindent {\small \sc Department of Mathematics, Rutgers University,
Piscataway, NJ 08854}

\noindent {\em E-mail address}: lepowsky@math.rutgers.edu

\vspace{1em}

\noindent {\small \sc Department of Mathematics, Rutgers University,
Piscataway, NJ 08854}

\noindent {\em E-mail address}: linzhang@math.rutgers.edu


\begin{thebibliography}{FGST2}

\bibitem[FHL]{FHL}
I.~B. Frenkel, Y.-Z. Huang and J.~Lepowsky,
On axiomatic approaches to vertex operator algebras and modules,
preprint, 1989;
{\em Memoirs Amer. Math. Soc.} {\bf 104}, 1993.

\bibitem[H]{tensor4}
Y.-Z. Huang, A theory of tensor products for module categories for a
vertex operator algebra, IV, {\em J. Pure Appl. Alg.} 100 (1995)
173--216.

\bibitem[HL1]{tensor1}
Y.-Z. Huang and J. Lepowsky, A theory of tensor products for module
categories for a vertex operator algebra, I, {\em Selecta Mathematica
(New Series)} {\bf 1} (1995), 699--756.

\bibitem[HL2]{tensor2}
Y.-Z. Huang and J. Lepowsky, A theory of tensor products for module
categories for a vertex operator algebra, II, {\em Selecta Mathematica
(New Series)} {\bf 1} (1995), 757--786.

\bibitem[HL3]{tensor3}
Y.-Z. Huang and J. Lepowsky, A theory of tensor
products for module categories for a vertex operator algebra, III,
{\em J. Pure Appl. Alg.} {\bf 100} (1995) 141--171.

\bibitem[HLZ1]{HLZ1} Y.-Z.~Huang, J.~Lepowsky and L.~Zhang, Logarithmic
tensor category theory for generalized modules for a conformal
vertex algebra, I: Introduction and strongly graded
algebras and their generalized modules, to appear.

\bibitem[HLZ2]{HLZ2} Y.-Z.~Huang, J.~Lepowsky and L.~Zhang, Logarithmic
tensor category theory, II: Logarithmic formal calculus
and properties of logarithmic intertwining operators, to appear.

\bibitem[HLZ3]{HLZ4} Y.-Z.~Huang, J.~Lepowsky and L.~Zhang, Logarithmic
tensor category theory, IV: Constructions of tensor
product bifunctors and the compatibility conditions, to appear.

\bibitem[HLZ4]{HLZ5} Y.-Z.~Huang, J.~Lepowsky and L.~Zhang, Logarithmic
tensor category theory, V: Convergence condition for
intertwining maps and the corresponding compatibility
condition, to appear.

\bibitem[HLZ5]{HLZ6} Y.-Z.~Huang, J.~Lepowsky and L.~Zhang, Logarithmic
tensor category theory, VI: Expansion condition, associativity of logarithmic
intertwining operators, and the associativity isomorphisms, to appear.

\bibitem[HLZ6]{HLZ7} Y.-Z.~Huang, J.~Lepowsky and L.~Zhang, Logarithmic
tensor category theory, VII: Convergence and extension
properties and applications to expansion for intertwining
maps, to appear.

\bibitem[HLZ7]{HLZ8} Y.-Z.~Huang, J.~Lepowsky and L.~Zhang, Logarithmic
tensor category theory, VIII: Braided tensor category
structure on categories of generalized modules for a
conformal vertex algebra, to appear.

\end{thebibliography}
\end{document}